\documentclass[5p]{elsarticle}

\usepackage{amsmath} 
\usepackage{amssymb}  
\usepackage{bm}  

\usepackage{xcolor}
\usepackage{tikz}
\usepackage{algorithm}
\usepackage{caption}  
\usepackage{booktabs}  
\usepackage{hyperref}
\usepackage[capitalise]{cleveref}  
\usepackage{url}  
\usepackage{mathtools}  

\graphicspath{{examples/}}
\makeatletter
\newcommand\appendtographicspath[1]{%
  \g@addto@macro\Ginput@path{#1}%
}
\makeatother

\newdefinition{remark}{Remark}
\newtheorem{theorem}{Theorem}
\newtheorem{proposition}{Proposition}
\newproof{pf}{Proof}

\newcommand{\cC}{\mathcal{C}}

\newcommand{\cCH}[1]{\mathcal{H}}
\newcommand{\cCV}[1]{\mathcal{V}}
\newcommand{\psiV}{\psi_\mathcal{V}}
\newcommand{\DpsiV}[2]{\psiV'(#1;#2)}

\newcommand{\cA}{\mathcal{A}}
\newcommand{\cM}{\mathcal{M}}

\newcommand{\innerproduct}[2]{#1^T #2}

\newcommand{\cG}{\mathcal{G}}
\newcommand{\R}{\mathbb{R}}

\newcommand{\Lone}{\ensuremath{\mathcal{L}_1}}

\newcommand{\Linfty}{\ensuremath{\mathcal{L}_\infty}}

\newcommand{\posvector}{\ensuremath{p}}

\newcommand{\rate}{\eta}
\newcommand{\numrays}{m}

\DeclareMathOperator{\interior}{int}

\DeclareMathOperator{\convexhull}{\mathbf{Co}}

\DeclareMathOperator{\offdiag}{off-diag}

\DeclareMathOperator{\vect}{vec}

\DeclareMathOperator{\sgn}{sgn}
\DeclareMathOperator{\argmin}{argmin}
\DeclareMathOperator{\argmax}{argmax}

\begin{document}

\title{
    Polyhedral Estimation of
    $\Lone$ and $\Linfty$ Incremental Gains of Nonlinear Systems
    \tnoteref{t1}
}
\tnotetext[t1]{This work was supported by the Engineering and Physical Sciences Research Council (EPSRC) of the United Kingdom - Project Reference 1950140.}

\author[1]{Dimitris Kousoulidis}
\ead{dk483@cam.ac.uk}

\author[1]{Fulvio Forni\corref{cor1}}
\ead{f.forni@eng.cam.ac.uk}
\cortext[cor1]{Corresponding author}

\address[1]{University of Cambridge, Department of Engineering, Trumpington Street, CB2 1PZ Cambridge, United Kingdom}

\begin{abstract}
We provide novel dissipativity conditions for bounding
the incremental $\Lone$ gain of systems.
Moreover, we adapt existing results on the $\Linfty$ gain
to the incremental setting and 
relate the incremental $\Lone$ and $\Linfty$ 
gain bounds through system adjoints.
Building on work on optimization based approaches to
constructing polyhedral Lyapunov functions,
we make use of these conditions
to obtain a Linear Programming based algorithm
that can provide increasingly sharp bounds
on the gains as a function of a given candidate
polyhedral storage function or polyhedral set.
The algorithm is also extended to
allow for the design of linear feedback controllers
for performance,
as measured by the bounds on the incremental gains.
We apply the algorithm to a couple of numerical examples
to illustrate the power, as well as some limitations,
of this approach.
\end{abstract}


\maketitle

\section{Introduction}
\label{sec:intro}

When we treat open dynamical systems as operators $\cG$
that map input signals $w$ onto output signals $z = \cG w$,
the induced operator norms of $\cG$
quantitatively measure the amplification factor of
the input-output relationship,
which we call the system gain.
Gains provide a useful tool for system specifications,
allowing us to quantitatively characterize performance
and to pose control design problems
to optimize them in closed loop
\cite{zhou_robust_1995}.
Gains are also used
for studying the stability of interconnections
through the small gain theorem,
which states that
the feedback interconnection of two systems $\cG_1$ and $\cG_2$ 
will be stable as long as $\|\cG_1\|\|\cG_2\| < 1$
\cite[Chapter 3]{desoer_feedback_1975}.
For these purposes,
gains based on various norms have been
extensively studied and used in systems and control theory.
The selection of the norm needs to balance the needs of capturing performance
versus our ability to compute with it
\cite{boyd_linear_1991,zhou_robust_1995}.
In this paper we will be using incremental gains.
Incremental gains are a classic concept in control literature \cite{desoer_feedback_1975}.
They bound the incremental amplification between any pair of input signals.
Incremental gains maintain the interpretations and properties of gains discussed above,
including a form of the small gain theorem
\cite[Chapter 3]{desoer_feedback_1975}.
For linear systems, gains and incremental gains are equivalent.

We will be computing bounds on
the incremental $\Lone$ and the $\Linfty$ gains:
\begin{equation}
\|\cG\|_{1} = \sup_{w_1-w_2 \in \Lone \setminus \{0\}} \frac{\|z_1-z_2\|_1}{\|w_1-w_2\|_1}
,
\label{eq:lone-incr-gain}
\end{equation}
and
\begin{equation}
\|\cG\|_{\infty} = \sup_{w_1-w_2 \in \Linfty \setminus \{0\}} \frac{\|z_1-z_2\|_\infty}{\|w_1-w_2\|_\infty}
,
\label{eq:linfty-incr-gain}
\end{equation}
respectively, where 
\begin{alignat*}{2}
&\|w\|_1 &&= \int_0^\infty |w|_1 dt = \int_0^\infty \sum_i |w_i| dt,
\\
&\|w\|_\infty &&= \sup_t |w|_\infty = \sup_{i,\,t} |w_i|
.
\end{alignat*}
$w \in \Lone$ if $\|w\|_1 < \infty$,
and $w \in \Linfty$ if $\|w\|_\infty < \infty$.

The $\Lone$ norm assigns lower magnitude to short, high amplitude signals
over long, low amplitude signals;
while the $\Linfty$ norm assigns lower magnitude to long, low amplitude signals
over short, high amplitude signals.
The $\Lone$ gain is most suitable in quantifying specifications
in which the signals are each directly related to the consumption of a resource,
for example fuel use in a rocket or transaction fees in a trading system,
while the $\Linfty$ gain is most suitable in quantifying specifications
about the worse-case amplitudes of signals.
Moreover, for linear systems,
upper bounds on the $\Lone$ and $\Linfty$ gains, $\gamma_1$ and $\gamma_\infty$,
provide upper bounds for any $\mathcal{L}_p$ gain with $1 \leq p \leq \infty$, $\gamma_p$
\cite{desoer_feedback_1975}
\begin{equation}
\gamma_p \leq \gamma_1^{1/p}\gamma_\infty^{1/q}
,\ \text{where } 1/p+1/q = 1
,
\label{eq:desoer-lp-bound}
\end{equation}
offering another reason to study this combination of gains.

$\Lone$ and $\Linfty$ gains have been mostly studied
in the context of positive systems.
This is
because of the ease with which they can be expressed and
their efficient calculation \cite{briat_robust_2013,rantzer_scalable_2021}.
This has also led to conditions for 
robust stability and stabilization via integral linear constraints,
a linear analogue to integral quadratic constraints
\cite{briat_robust_2013}.
Moreover,
weighted $\Lone$ gains can be used to obtain a complete characterization of stability
for positive interconnections of positive systems \cite{ebihara_l1_2011}.

Beyond positive systems,
$\Lone$ and $\Linfty$ gains of linear systems are computed
through numerical integration of the impulse response
\cite{balakrishnan_computing_1992,linnemann_computing_2007,rutland_computing_1995}.
Current results for $\Lone$ and $\Linfty$ bounds on gains
for nonlinear systems include
\cite[Theorem 5.1]{khalil_nonlinear_2002},
where a uniform bound for all $\mathcal{L}_p$ gains
in terms of a provided Lyapunov function
is derived.
The work of \cite{blanchini_set-theoretic_2015}
and \cite{fialho_l1_1995}
presents algorithms for computing the $\Linfty$ gain of
linear time/parameter varying (LTV/LPV) systems
via invariant polyhedral sets.

The main contributions of this paper are novel conditions
for computing incremental $\Lone$ gain bounds based
on dissipation inequalities with respect to polyhedral storage functions.
Conditions for incremental $\Linfty$ gain bounds
using polyhedral set invariance
are also presented.
The $\Linfty$ conditions are closely related to prior work
\cite{blanchini_set-theoretic_2015,fialho_l1_1995},
but tweaked for the incremental gain setting.
Both sets of conditions are framed as optimization problems
and solutions to either one of the gains
can be mapped to solutions of the other for an adjoint system,
allowing for any algorithm designed
to compute incremental $\Lone$ gains to be used
for computing incremental $\Linfty$ gain bounds and vice versa.
This generalizes the known relationship between
$\Lone$ and $\Linfty$ gains for linear systems
\cite{briat_robust_2011}.
When the polyhedral function/set is fixed
the conditions are reduced to linear programming problems.
This allows us to adapt our previous work
on finding polyhedral Lyapunov functions of fixed complexity
\cite{kousoulidis_polyhedral_2021}
to this problem
for both analysis and synthesis.
We provide definitions and the conditions at an abstract level
in \cref{sec:gains},
specialize them to polyhedral functions and provide our novel
conditions in \cref{sec:poly},
provide an overview of our adapted analysis and synthesis
algorithms in \cref{sec:alg},
and apply the algorithms on some numerical examples
in \cref{sec:examples}.

{
\small
\vspace{3mm}
\textbf{Notation:}
We use calligraphic letters to denote sets and operators,
and capital letters to denote matrices.
We use $\bm{1}$ to denote vectors of all ones and
$[k]$ to denote set $\{1,\dots,k\}$.
Operations applied to a set are meant to denote operations
applied to each element of the set;
for example,
we use $\cA^T$ to denote $\{A^T : A \in \cA\}$.
$A^i$ denotes the $i^{th}$ row of matrix $A$,
while
$A_j$ denotes the $j^{th}$ column.
$[A,B]$ denotes the matrix obtained by stacking matrices $A$ and $B$ \emph{horizontally},
while
$[A;B]$ denotes the matrix obtained by stacking matrices $A$ and $B$ \emph{vertically}.
The vectorization operator $\vect(\cdot)$
converts
a $n \times m$ matrix into
a column vector with $nm$ entries
by stacking up its columns.
We use $\|x\|$ to denote the norm of a \emph{signal} $x$
and $|x|$ to denote the norm of a \emph{vector} $x$.
The sign function
$\sgn(x)$ is set to $-1 \text{ if } x < 0, 0 \text{ if } x = 0, \text{and } 1 \text{ if } x > 0$.
For any given matrix $A$, $A \geq 0$ means that all the elements of $A$ are non-negative.
}
\section{Gains and Incremental Gains}
\label{sec:gains}

\subsection{Gains for Linear Systems}
\label{sec:gains_linear}
For linear systems,
the $\Lone$ and $\Linfty$ gains
can be characterized in terms of
the matrix of impulse responses $H$.
We denote 
the entry in the $i^{th}$ row and $j^{th}$ column of $H$
as $h_{ij}$.
This corresponds to the impulse response of
the $i^{th}$ output to the $j^{th}$ input.

The $\Lone$ and $\Linfty$ gains of a system with $m$ outputs and $n$ inputs are then
\cite{briat_robust_2013}
\begin{equation}
   \|\cG\|_1 = \max_{1 \leq j \leq n} \left( \sum_{i = 1}^{m} \left( \int_0^\infty |h_{ij}(t)| dt \right) \right)
   \label{eq:lone-tf}
\end{equation}
and
\begin{equation}
   \|\cG\|_\infty = \max_{1 \leq i \leq m} \left( \sum_{j = 1}^{n} \left( \int_0^\infty |h_{ij}(t)| dt \right) \right)
   .
   \label{eq:linfty-tf}
\end{equation}
Note that,
other than for the computation of the $m \times n$ integrals $\int_0^\infty |h_{ij}(t)| dt$,
the remaining computations for \eqref{eq:lone-tf} and \eqref{eq:linfty-tf}
amount to summations and maximizations over a finite number of scalars.
As such evaluating $\int_0^\infty |h_{ij}(t)| dt$ numerically is the main computational challenge
for evaluating the two gains
and various methods have been proposed for computing it
\cite{balakrishnan_computing_1992,linnemann_computing_2007,rutland_computing_1995}.

It follows from \eqref{eq:lone-tf} and \eqref{eq:linfty-tf} that
the $\Lone$ gain of a system with impulse responses $H$
is equivalent to the $\Linfty$ gain of a system with impulse responses $H^T$.
This also means that 
the $\Lone$ and $\Linfty$ gains  
of single-input single-output or self-adjoint linear systems
are equivalent.

\begin{remark}
    \label{remark:pos}
For the special case of positive linear systems,
all $h_{ij}(t) \geq 0$ for all $t \geq 0$
and so $\int_0^\infty |h_{ij}(t)| dt = \int_0^\infty h_{ij}(t) dt$,
which is also equivalent to the frequency response of $\cG$ evaluated at $0$.
This enables efficient calculation and design for
$\Lone$ and $\Linfty$ gains \cite{briat_robust_2013,rantzer_scalable_2021}.
\end{remark}

\subsection{Incremental Gains for Nonlinear Systems}
The approach in \cref{sec:gains_linear} relies
on concepts from linear time invariant systems theory that are hard to generalize
to other classes of systems.
A more general approach to computing incremental gains
(or, at least, upper bounds on incremental gains)
is through dissipativity theory
\cite{willems_dissipative_1972},
which we apply to more general,
minimal realization
state-space systems of the form
\begin{equation}
\dot{x} = f(x) + Bw \qquad z = Cx
,
\label{eq:open-state-space}
\end{equation} 
with $n$ states, $n_w$ disturbance inputs, and $n_z$ outputs.
To obtain a well defined operator we also set $x(0) = 0$.

Under this framework we replace conditions over all signals,
such as those in \eqref{eq:lone-incr-gain} and \eqref{eq:linfty-incr-gain},
with point-wise in time conditions over all time, inputs, and state variables
-- called \emph{dissipation inequalities}.
This framework then requires finding a positive function of the state,
the \emph{storage function},
which plays a central role in the dissipation inequalities.
For simplicity,
in this paper we only consider linear input/output system matrices $B$ and $C$.
However, the approach can be in principle be further generalized
to nonlinear dependencies as well.

For incremental $\Lone$ gains
the suitable dissipation inequality for all
$w_1,\, w_2$ and $x_1,\, x_2$
is
\begin{equation}
    \dot{\psi} \leq - |z_1-z_2|_1 + \gamma_1 |w_1-w_2|_1 
    ,
    \label{eq:lone-dissipation-gen}
\end{equation}
where $\psi(x_1,x_2)$ is a storage function.
To make it simpler to verify \eqref{eq:lone-dissipation-gen},
we limit ourselves to considering storage functions of the form $\psi(x_1-x_2)$
that satisfy $\psi(x_1-x_2) \geq 0$ for all $x_1, x_2 \in \mathbb{R}^n$, and $\psi(0) = 0$.
Integrating \eqref{eq:lone-dissipation-gen}
in time from $0$ to $T$
gives
\begin{equation*}
    \|z_1-z_2\|_1 \leq -\psi(x_1(T)-x_2(T)) + \gamma_1 \|w_1-w_2\|_1
\end{equation*}
and so, since $\psi(x_1(T)-x_2(T)) \geq 0$,
\begin{equation*}
    \frac{\|z_1-z_2\|_1}{\|w_1-w_2\|_1} \leq \gamma_1
    ,
\end{equation*}
giving us $\gamma_1$ as an upper bound for the incremental $\Lone$ gain of the system.

For $\Linfty$ there is no need to integrate the signals involved
because they directly correspond to point-wise in time conditions
which are conveniently expressed in terms of set invariance
\cite[Section 6.4]{blanchini_set-theoretic_2015}.
What we need is a limit on the maximum incremental amplification
of the output relative to the input.
We capture this 
by confirming that inputs with bounded incremental norms lead to
states with bounded difference,
i.e. confined to a bounded set $\mathcal{X}$,
for all $t \geq 0$,
\begin{subequations}
\begin{align}
\|w_1-w_2\|_\infty \leq 1
,\,
x_1(0) &= 0
,\,
x_2(0) = 0
\nonumber
\\
&\implies (x_1(t)-x_2(t)) \in \mathcal{X}
\label{eq:linfty-input-state}
\\
(x_1(t)-x_2(t)) \in \mathcal{X} &\implies \|z_1-z_2\|_\infty \leq \gamma_\infty
,
\label{eq:linfty-state-output}
\end{align}
\label{eq:linfty-gen}
\end{subequations}
giving $\gamma_\infty$
as an upper bound for the incremental $\Linfty$ gain of the system.
Here the set $\mathcal{X}$ plays a role similar to that of a storage function.

We next specialize
\eqref{eq:lone-dissipation-gen}
and \eqref{eq:linfty-gen}
to the family of polyhedral functions and sets respectively.

\section{Polyhedral Approximation of Gains}
\label{sec:poly}

\subsection{Polyhedral Preliminaries}
To compute gains we use C-sets and gauge functions.
C-sets are compact convex sets that contain $0$ in their interior,
and gauge functions are subadditive positively homogeneous functions (of degree $1$).
We can extract a gauge function from a C-set through its Minkowski functional and
a C-set from a gauge function through its unit ball,
so each C-set induces a gauge function and vice versa \cite{blanchini_set-theoretic_2015}.

Polyhedra have 
two distinct representations,
which we call V-representation and H-representation,
V-rep and H-rep for short.
A V-rep polyhedral set is represented by its extremal vertices,
while a H-rep polyhedral set
is represented by intersections of affine half-spaces.
We adopt the convention of using
$\cCV{V}$ to denote V-rep sets
constructed from numerical matrix $V$.
Likewise,
we use $\cCH{H}$ to denote H-rep sets
constructed from numerical matrix $H$.
\begin{description}
    \item[V-rep] given $n \times \numrays$ matrix $V$,
    \begin{equation}
        \cCV{V}
        = \{x \in \R^n: x = V\posvector,\, \posvector \geq 0,\, \innerproduct{\bm{1}}{\posvector} = 1\}
        \label{eq:V-rep}
    \end{equation}
    with interior
    \begin{equation}
        \interior \cCV{V} = \{x: x = V\posvector,\, \posvector > 0,\, \innerproduct{\bm{1}}{\posvector} = 1\}
        \label{eq:V-rep_int}
    \end{equation}
    and corresponding gauge function
    \begin{subequations}
        \label{eq:psiV}
    \begin{align}
        \psiV(x)
        = &\min_{p \in \R^m}\{\bm{1}^T p : x = Vp,\, p \geq 0 \}
        \label{eq:psiV_primal} \\
        = &\max_{h \in \R^n}\{h^T x : h^T V \leq \bm{1}^T \}
        \label{eq:psiV_dual} \,
    \end{align}
    \end{subequations}
    where \eqref{eq:psiV_dual} is the dual problem to \eqref{eq:psiV_primal}
    \cite{boyd_convex_2004}.
    \item[H-rep] given $\numrays \times n$ matrix $H$,
    \begin{equation}
        \cCH{H} = \{x \in \R^n: Hx \leq \bm{1} \}
        \label{eq:H-rep}
    \end{equation}
    with interior
    \begin{equation}
        \interior \cCH{H} = \{x: Hx < \bm{1} \}
        .
        \label{eq:H-rep_int}
    \end{equation}
\end{description}

$\cCV{V}$ is the 
convex hull of the columns of $V$,
which are called the vertices of $\cCV{V}$.
On the other hand, $\cCH{H}$ is the intersection of half-spaces based on the
rows of $H$.
We use gauge functions of V-rep sets for computing $\Lone$ gains
and H-rep sets for computing $\Linfty$ gains.

We take a set-theoretic approach to computing 
$\Linfty$ gains
so we omit the functional perspective for brevity.
On the other hand,
because we compute the $\Lone$ gain through dissipation inequalities,
we need to characterize the rate of change of V-rep gauge functions,
which will act as the candidate storage functions in \eqref{eq:lone-dissipation-gen}.
We do this through the use of subdifferentials \cite{rockafellar_convex_1970}:
\begin{equation}
\partial\psi(x) = \{h : \forall v\in \R^n, \, \psi(v)-\psi(x) \geq h^T(v-x) \}
.
\label{eq:subdiff}
\end{equation}
To provide some intuition about subdifferentials,
the subdifferential of an everywhere differentiable convex function
consists at every point of a unique element
that is equivalent to its gradient.
Unfortunately,
$\psiV(x)$ is not everywhere differentiable
but we can still obtain a convenient numerical representation of $\partial\psiV(x)$.
\begin{proposition}[Subdifferential of $\psiV$]
    \label{prop:subdiff-v}
    \begin{equation}
        \partial\psiV(x) = \{h : h^T x = \psiV(x),\,h^TV \leq \bm{1}^T \}
        .
        \label{eq:subdiff-v}
    \end{equation}
\end{proposition}
\begin{pf}
We first rewrite $\psiV$ in the following equivalent form
$$
    \psiV(x) = \max_{h \in \mathcal{Z}} h^T x
$$
where $\mathcal{Z} = \{h : V^T h \leq \bm{1}\}$. We then apply Danskin’s Theorem
  \cite[Proposition 4.5.1]{bertsekas_convex_2003} to derive the expression of the subdifferential.
  To apply the theorem
  we must show that
  $\mathcal{Z}$ is compact.
  $\mathcal{Z}$ is the intersection
  of closed half-spaces so is closed and,
  using \cite[Corollary 14.5.1]{rockafellar_convex_1970},
  it must also be bounded
  because
  $0 \in \interior \cCV{V}$
  and $\cCV{V}$ is a C-set.
  From \cite[Proposition 4.5.1]{bertsekas_convex_2003}, the subdifferential is then given by the convex hull
  of all $h$ such that $h \in \mathcal{Z}$ and maximize $h^T x$. The latter condition corresponds to $h^Tx = \psiV(x)$.
  Combining these constraints together we get the set
  $\{ h : h^Tx = \psiV(x),\, h^T V \leq \bm{1}^T\}$, which is already convex. This leads to \eqref{eq:subdiff-v}.
\hfill $\Box$
\end{pf}
Given any direction $z \in \R^n$,
$\psiV(x)$ also has a well defined directional derivative
\cite[Theorem D1.2.2]{hiriart-urruty_fundamentals_2001}
\begin{subequations}
\begin{align}
\DpsiV{x}{z}
\!&= \max_{h \in \partial\psiV(x)}\ h^T z
\nonumber
\\
&=\max_{h}\{h^T z : h^T x = \psiV(x),\,h^TV \leq \bm{1}^T \}
\label{eq:derivative_v}
\\
&=\min_{p,k}\{\bm{1}^Tp \!+\! \psiV(x)k : z \!=\! Vp \!+\! kx, p \!\geq\! 0\}
\label{eq:derivative_v_dual}
\end{align}
\label{eq:derivative_v_both}
\end{subequations}
where \eqref{eq:derivative_v_dual} is the dual problem to \eqref{eq:derivative_v}.

We can integrate \eqref{eq:derivative_v_both}
to obtain the total change in $\psiV$
under some dynamics $\dot{x} = f(x)$
\begin{equation}
    \psiV(x(t)) - \psiV(x(0)) = \int_0^t
    \DpsiV{x(\tau)}{f(x(\tau))} d\tau
    .
\label{eq:integral_general}
\end{equation}
This follows from
\cite[Theorem 2.10 and 2.11]{blanchini_set-theoretic_2015};
see also \cite[Appendix 1 Theorem 4.3]{rouche_stability_1977} and 
\cite[Corollary 24.2.1]{rockafellar_convex_1970}.

\cref{eq:derivative_v} shows that the rate of change
of V-rep polyhedral functions can be characterized
in terms of linear programming.

\subsection{Polyhedral Gain Conditions}
\label{sec:poly-conditions}

To compute bounds on incremental $\Lone$ and $\Linfty$ gains
we will use the variational system of \eqref{eq:open-state-space},
along generic solution pairs $(x,w):\mathbb{R} \to \mathbb{R}^n \times \mathbb{R}^{n_w}$ of \eqref{eq:open-state-space}:
\begin{equation}
    \dot{\delta x} = \partial f(x) \delta x + B \delta w ,\, \delta z = C \delta x .
    \label{eq:open-state-space-var}
\end{equation}
Along a generic solution $(x,w)$, \eqref{eq:open-state-space-var} 
corresponds to a time-varying linear system with $n$ (variational) states, $n_w$ inputs, and $n_z$ outputs, 
which captures the effect of (infinitesimal) perturbations on $(x,w)$ \cite{crouch_variational_1987,forni_differential_2014,forni_differential_2013}.

In what follows, we will take advantage of the fact that a uniform bound $\gamma$
that satisfies $\|\delta z\| \leq \gamma \|\delta w \|$ 
for all input/output pairs of the variational system
and uniformly with respect to all solution pairs $(x,w)$ of \eqref{eq:open-state-space}, 
is also a bound on the incremental gain of \eqref{eq:open-state-space}.
For instance, consider a generic pair of inputs signal $w_0$ and $w_1$, 
and define the parameterized signal $w_s = sw_1 + (1-s) w_0$, for $0\leq s\leq 1$.
Using \eqref{eq:open-state-space}, compute the state trajectory $x_s$ and the output $z_s$ 
related to the input $w_s$ and to the initial condition $x_s(0)= 0$. 
For any $0 \leq s \leq 1$, the quantities $\frac{d}{ds} w_s$ and $\frac{d}{ds} z_s$ correspond to the input and output signals of 
\eqref{eq:open-state-space-var} computed along the solution pair $(x_s,w_s)$. Thus,
\begin{equation}
\begin{split}
\|z_1 - z_0 \|  &= \left\| \int_0^1  \frac{d}{ds} z_s ds \right\|  \ 
\leq \int_0^1  \left\| \frac{d}{ds} z_s \right\| ds \\
& \leq \gamma  \int_0^1 \left\| \frac{d}{ds} w_s \right\| ds  \ 
= \gamma  \int_0^1 \| w_1 - w_0 \| ds  \\
&= \gamma  \| w_1 - w_0 \| , 
\end{split} 
\end{equation}
where the second inequalities follows from the uniform bound $\|\delta z\| \leq \gamma \|\delta w \|$.

Together with the use of the variational system,
we assume the
existence of a V-rep polyhedral relaxation of its differential dynamics
\begin{equation}
    \label{eq:differential-polyhedral-relaxation}
\cA = \{A_1,\dots,A_k\}
\ \text{such that, for all $x$,} \
    \partial f(x) \in \convexhull \cA
    ,
\end{equation}
where
\begin{equation*}
    \convexhull \cA
    = \left\{ \sum_{i=1}^k \lambda_i A_i : A_i \in \cA,\, \lambda_i \geq 0,\, \sum_{i=1}^k \lambda_i = 1 \right\}.
\end{equation*}
Namely, this means that, for any $x$, $\partial f(x) = \sum_{i=1}^k \lambda_i(x) A_i$ where $A_i \in \mathcal{A}$, $\lambda_i(x) \geq 0$, and 
$\sum_{i=1}^k \lambda_i(x) = 1$.
This relaxation is inspired by a standard step in the literature
for computation of nonlinear gains by
linear matrix inequalities (LMIs) \cite{boyd_linear_1994}.
We don't present a general algorithm for computing such relaxations
since it is a non-trivial problem.
This is related to the trade-off between the
quality of the relaxation and the number of elements in it.
Instead, in our examples,
we work directly with sets $\cA$,
with the results then applying to any nonlinear system
that satisfies \eqref{eq:differential-polyhedral-relaxation}.

For a linear system,
\eqref{eq:open-state-space-var} has the same dynamics
as \eqref{eq:open-state-space} with $f(x) = \partial f(x) = A$
and so $\cA = \{A\}$.

\begin{theorem}[Bound on incremental $\Lone$ gain]\ \\
\label{thm:incr-lone-gain}
Consider variables
$
\rate_w \in \R,\,
\rate_z \in \R,\,
V \in \R^{n \times m},\,
P \in \R^{m \times 2n_w},\,
\cM = \{M_1 \in \R^{m \times m},\dots,M_k \in \R^{m \times m}\}
$
(where $k$ is the size of $\cA$)
and vector valued function of $V$,
$$\hat{z}_V = [|CV_1|_1,\dots,|CV_m|_1]$$
(where $V_j$ is the $j^{th}$ column of $V$).

An upper bound for the incremental $\Lone$ gain of (8) is given by
    \begin{subequations}
        \label{eq:incr-lone-hard}
        \begin{align}
        \min_{\rate_w, \rate_z, V, P, \cM}
        \quad &\> \rate_w/\rate_z
        \nonumber
        \\
        \noalign{\medskip}
        \begin{split}
            P \geq 0,
            \    
            [B, -B] = VP,
            &\ 
            \rate_w \bm{1}^T = \bm{1}^T P
        \end{split}
        \label{eq:incr-lone-hard-input}
            \\
        \noalign{\medskip}
            \rate_z > 0,
            \
            \offdiag(M_i) \geq 0,
            \ 
            A_i &V = V M_i ,
            \    
            -\rate_z \hat{z}_V = \bm{1}^T M_i ,
            \nonumber
            \\
            \text{for all } &i \in [k]
            .
        \label{eq:incr-lone-hard-state}
        \end{align}
    \end{subequations}
    Moreover, the
    gauge function
    $\psi = \psiV/\rate_z$ satisfies \eqref{eq:lone-dissipation-gen}
    for $\gamma_1 = \rate_w/\rate_z$.
    \hfill $\lrcorner$
\end{theorem}

\begin{theorem}[Bound on incremental $\Linfty$ gain] \ \\
    \label{thm:incr-linfty-gain}
Consider variables
$
\rate_z \in \R,\,
\rate_w \in \R,\,
H \in \R^{m \times n},\,
P \in \R^{2n_z \times m},\,
\cM = \{M_1 \in \R^{m \times m},\dots,M_k \in \R^{m \times m}\}
$
and vector valued function of $H$,
$$\hat{w}_H = [|B^T(H^1)^T|_1;\dots;|B^T(H^m)^T|_1]$$
(where $H^j$ is the $j^{th}$ row of $H$).

An upper bound for the incremental $\Linfty$ gain of (8) is given by
    \begin{subequations}
        \label{eq:incr-linfty-hard}
        \begin{align}
        \min_{\rate_z, \rate_w, H, P, \cM}
        \quad &\> \rate_z/\rate_w
        \nonumber
        \\
        \noalign{\medskip}
            \rate_w > 0,
            \
            \offdiag(M_i) \geq 0,
            \    
            H &A_i = M_i H,
            \ 
            -\!\rate_w \hat{w}_H = M_i \bm{1},
            \nonumber
            \\
            \text{for all } &i \in [k]
        \label{eq:incr-linfty-hard-state}
            \\
        \noalign{\medskip}
        \begin{split}
            P \geq 0,
            \    
            [C; -C] = PH,
            &\ 
            \rate_z \bm{1} = P \bm{1}.
        \end{split}
        \label{eq:incr-linfty-hard-output}
        \end{align}
    \end{subequations}
    Moreover, the
    polyhedral set
    $\mathcal{X} = \cCH{H}/\rate_w$
    satisfies \eqref{eq:linfty-gen}
    for $\gamma_\infty = \rate_z/\rate_w$.
    \hfill $\lrcorner$
\end{theorem}

The columns of the matrix $V$ in \eqref{eq:incr-lone-hard} are the vertices of the 
level set of the storage function. The intuition is that \eqref{eq:incr-lone-hard-input} quantifies the impact that the 
(incremental) input has on the growth of the storage function. In the proof of Theorem \ref{thm:incr-lone-gain} this is well
captured by \eqref{eq:incr-lone-dissipation-u-only}. Likewise, \eqref{eq:incr-lone-hard-state} characterizes how the internal dynamics
of the system and its output features affect the storage, as formalized below by \eqref{eq:incr-lone-dissipation-x-only}. This splitting
is made possible by the use of storages with polyhedral representation. In a similar way, the rows of the matrix $H$ in \eqref{eq:incr-linfty-hard}
characterize a polyhedral set on the system state space as an intersection of half-spaces.
The goal of \eqref{eq:incr-linfty-hard-state} is to bound the expansion of this polyhedral set under the action of the 
forced system dynamics, in the sense of \eqref{eq:linfty-input-state}. Later, \eqref{eq:incr-linfty-hard-output} maps this bound
on the system state into a bound on the output, as clarified by \eqref{eq:linfty-state-output}. 
We refer the reader to the proofs of the theorems for a more detailed discussion on the connection between the conditions 
of the theorems and the incremental gains of the system.

The optimization problems in
\eqref{eq:incr-lone-hard} and \eqref{eq:incr-linfty-hard} are adjoint.
We mean this in the sense that
each solution of \eqref{eq:incr-lone-hard}
for system matrices $\{\cA,B,C\}$
with variables $(\rate_w, \rate_z, V, P, \cM)$
is also a solution of \eqref{eq:incr-linfty-hard}
for system matrices $\{\cA^T,C^T,B^T\}$ (the adjoint system)
with variables $(\rate_z, \rate_w, V^T, P^T, \cM^T)$
and vice versa.
This matches the expected relationship between the two gains
for linear systems 
(\eqref{eq:lone-tf} and \eqref{eq:linfty-tf}).
It also means that 
it is sufficient to implement only one of
\eqref{eq:incr-lone-hard} or \eqref{eq:incr-linfty-hard}.

\begin{remark}
Because the approach relies on the set inclusion of the variational dynamics \eqref{eq:open-state-space-var} into a polytopic set $\convexhull \cA$,
the same methodology can also be directly applied to
any system whose dynamics can be characterized by
polytopic linear differential inclusions (LDIs),
including linear time and parameter varying (LTV/LPV) systems.
\end{remark}

\begin{remark}
Both 
\eqref{eq:incr-lone-hard} and \eqref{eq:incr-linfty-hard}
make use of two positive constants instead of one,
which might appear unexpected.
This is because, in each case, one of the two constants
($\rate_w$ in \eqref{eq:incr-lone-hard-input} and
$\rate_z$ in \eqref{eq:incr-linfty-hard-output})
act as normalizers for the scale of the storage functions/sets.
The normalizing constants can be fixed to $1$ with no loss 
of generality
but we keep the additional degrees of freedom to better match
the presentation in \cref{sec:alg}.
\end{remark}

\subsection{Proofs}

\noindent \emph{Proof of \cref{thm:incr-lone-gain}:}
Let $w = w_1 - w_2$,
$x = x_1 - x_2$,
and $z = z_1 - z_2$.
    Our aim is to show that every solution $(\rate_w, \rate_z, V, P, \cM)$
    of \eqref{eq:incr-lone-hard}
    satisfies,
    for all $w$ and $x$
    \begin{equation}
        \dot{\psi}_{\cCV{V}} \leq -\rate_z|z|_1 + \rate_w |w|_1
        \label{eq:incr-lone-dissipation}
    \end{equation}
    for system \eqref{eq:open-state-space} and
    subject to \eqref{eq:differential-polyhedral-relaxation},
    where
    $$
    \dot{\psi}_{\cCV{V}} = \DpsiV{x}{f(x_1)-f(x_2)+Bw}
    .
    $$
    But, from \eqref{eq:derivative_v},
    \begin{align*}
    \DpsiV{x}{f(x_1)-&f(x_2) + B w} \\
    \leq \DpsiV{x}{&f(x_1)-f(x_2)} + \DpsiV{x}{B w}
    ,
    \end{align*}
    so \eqref{eq:incr-lone-dissipation} is implied by
    \begin{equation}
        (\forall w,x) \
        \DpsiV{x}{B w}
        - \rate_w |w|_1 \leq 0
        \label{eq:incr-lone-dissipation-u-only}
    \end{equation}
    and
    \begin{equation}
        (\forall x_1,x_2) \
        \DpsiV{x}{f(x_1)-f(x_2)} + \rate_z |C x |_1 \leq 0
        ;
        \label{eq:incr-lone-dissipation-x-only}
    \end{equation}
    which we address separately.

    For \eqref{eq:incr-lone-dissipation-u-only},
    we have that
    \begin{align*}
    \DpsiV{x}{B w} &=
    \max_{h}\{h^T B w : h^T x = \psiV(x),h^T V \leq \bm{1}^T \}
    \\
    &\leq \max_{h}\{h^T B w : h^T V \leq \bm{1}^T \}
    = \psiV(Bw)
    .
    \end{align*}
    Treating the 1-norm as the polyhedral gauge function
    with V-rep matrix
    $[I, -I]$,
    for any $w$ there exists a $p^* \geq 0$ such that
    $w = [I, -I]p^*$ and $\bm{1}^Tp^* = |w|_1$.
    Then, from \eqref{eq:incr-lone-hard-input},
    $Bw = [B, -B]p^* = VPp^*$, $(Pp^*) \geq 0$.
    Finally, from \eqref{eq:psiV_primal} and \eqref{eq:incr-lone-hard-input}
    \begin{equation*}
    \psiV(B w) \leq \bm{1}^TPp^* = \rate_w |w|_1
    .
    \end{equation*}

    To show \eqref{eq:incr-lone-dissipation-x-only}
    we first bound $\DpsiV{x}{f(x_1)-f(x_2)}$ in terms of $\cA$.
    Recall that $\partial\psiV(x) = \{ h : h^Tx = \psiV(x), h^TV \leq \bm{1}^T\}$.
    Then
    \begin{align*}
       &\DpsiV{x}{f(x_1)-f(x_2)} = \max_{h \in \partial\psiV(x)}\ h^T\big(f(x_1)-f(x_2)\big)
        \\
        &= \max_{h \in \partial\psiV(x)}\ h^T\left(\int_0^1 \partial f\big(sx_1 + (1-s)x_2\big)x ds
        \right)
        \\
        &\leq \int_0^1 \max_{h_s \in \partial\psiV(x)}\ h_s^T\big(\partial f(sx_1 +(1-s)x_2)x\big)ds
        \\
        &= \int_0^1 \max_{h_s \in \partial\psiV(x)}\ h_s^T\left(\sum_{i=1}^k\lambda_i(x_1,x_2,s)A_ix\right) ds
        \\
        &\leq \int_0^1 \sum_{i=1}^k \lambda_i(x_1,x_2,s) \max_{h_i \in \partial\psiV(x)}\ h_i^T\left(A_ix\right) ds
        \\
        &\leq \int_0^1 \left(\sum_{i=1}^k \lambda_i(x_1,x_2,s)\right)
        \max_{i} \left(\max_{h_i \in \partial\psiV(x)}\ h_i^TA_ix \right) ds
        \\
        &= \max_{i}\left(\max_{h_i \in \partial\psiV(x)}\ h_i^TA_ix\right)
        = \max_{i}\left(\DpsiV{x}{A_i x}\right)
        .
    \end{align*}
    From \eqref{eq:psiV_primal},
    for any $x$ there exists a $p^* \geq 0$ such that
    $x = Vp^*$ and $\bm{1}^Tp^* = \psiV(x)$.
    Combining with \eqref{eq:incr-lone-hard-state}
    and rewriting $M_i = Q_i + k^*I$ for some $Q_i \geq 0$ and some $k^*\in \mathbb{R}$,
    we have
    $A_i x = A_iVp^* = VM_i p^* = V(Q_i+k^*I)p^* = VQ_ip^* + k^*x$, where $(Q_ip^*) \geq 0$.
    Finally, from \eqref{eq:derivative_v_dual},
    \eqref{eq:incr-lone-hard-state},
    and the subadditivity of the 1-norm
    \begin{align*}
    \DpsiV{x}{A_i x} 
    &= \min_{p,k}\{\bm{1}^Tp \!+\! \psiV(x)k : A_i x \!=\! Vp \!+\! kx, p \!\geq\! 0\} \\
    &\leq \bm{1}^TQ_ip^* \!+\! \bm{1}^Tp^*k^*  
    = \bm{1}^T M_i p^* \\
    &= - \rate_z \sum_j |CV_j|_1 p_j^* 
    \leq - \rate_z |Cx|_1 .
    \end{align*} 
\hfill $\Box$

\noindent \emph{Proof of \cref{thm:incr-linfty-gain}:}
This is a straightforward variant of a known result
in the closely related LTV/LPV setting \cite{blanchini_set-theoretic_2015},
which we prove here for completeness.

Let $w = w_1 - w_2$,
$x = x_1 - x_2$,
and $z = z_1 - z_2$.
Our aim is to show that every solution $(\rate_z, \rate_w, H, P, \cM)$
of \eqref{eq:incr-linfty-hard} satisfies,
for all $t \geq 0$,
\begin{subequations}
\begin{align}
\|w\|_\infty \leq \rate_w, x(0) = 0
&\implies x(t) \in \cCH{H}
\label{eq:linfty-input-state-final}
\\
x(t) \in \cCH{H} &\implies \|z\|_\infty \leq \rate_z
.
\label{eq:linfty-state-output-final}
\end{align}
\label{eq:linfty-gen-final}
\end{subequations}
We first show that \eqref{eq:incr-linfty-hard-state} implies \eqref{eq:linfty-input-state-final}.
Then we show that \eqref{eq:incr-linfty-hard-output} implies \eqref{eq:linfty-state-output-final}.

To show \eqref{eq:linfty-input-state-final}
we write the conditions for the forward invariance of $\cCH{H}$
subject to the incremental dynamics
$\dot{x} = f(x_1)-f(x_2) + Bw$ and $|w|_\infty \leq \rate_w$.
Let
$
\cC_j = \{ x : H^jx = 1, Hx \leq \bm{1}\}
$.
Via Nagumo's Theorem \cite[Section 4.2]{blanchini_set-theoretic_2015},
\eqref{eq:linfty-input-state-final} holds if
\begin{equation}
\max_{x_1-x_2 \in \cC_j,|w|_\infty \leq \rate_w}\ H^j (f(x_1)-f(x_2)) + H^j B w
\leq 0
\label{eq:incr-linfty-state-invariance-definition}
\end{equation}
for all $j \in [m]$.
In fact, the meaning of the above inequality is that
at each point $x = x_1-x_2$ on the boundary of $\cCH{H}$,
the vector field $\dot{x} = f(x_1)-f(x_2) + Bw$ is directed towards the inside of $\cCH{H}$.
This guarantees that $\mathcal{H}$ is a forward invariant set for the incremental dynamics
for any $\|w\|_\infty \leq \rate_w$.

\cref{eq:incr-linfty-state-invariance-definition}
holds if the following inequality holds
\begin{equation}
\max_{x_1-x_2 \in \cC_j}\ H^j (f(x_1)-f(x_2))
\leq 
- \rate_w|B^T(H^j)^T|_1
, 
\label{eq:incr-linfty-state-invariance-inequality}
\end{equation}
since, from Hölder's inequality,
$$|H^j\! Bw| \leq \sum |(H^j\!B)^T_i w_i| \leq \rate_w |B^T\!(H^j)^T|_1 \ .$$

We now prove that \eqref{eq:incr-linfty-hard-state} implies \eqref{eq:incr-linfty-state-invariance-inequality}
(and thus \eqref{eq:linfty-input-state-final}). First,
\begin{align*}
&\max_{x_1-x_2 \in \cC_j}\ H^j (f(x_1)-f(x_2)) 
\nonumber
\\
&= \max_{x_1-x_2 \in \cC_j}\ 
H^j\left(\int_0^1 \partial f(sx_1 + (1-s)x_2)x ds\right)
    \nonumber
    \\
&= \max_{x_1-x_2 \in \cC_j}\ 
H^j\left(\int_0^1 \sum_{i=1}^k\lambda_i(x_1,x_2,s)A_ix ds\right)
    \nonumber
    \\
&= \max_{x_1-x_2 \in \cC_j}\ 
\left(\int_0^1 \sum_{i=1}^k\lambda_i(x_1,x_2,s)H^jA_ix ds\right)
    \nonumber
\end{align*}
\begin{align*}
&\leq \max_{x_1-x_2 \in \cC_j}\ 
\int_0^1 \left(\sum_{i=1}^k\lambda_i(x_1,x_2,s)\right)
\max_{i}\left(H^jA_ix\right) ds
    \nonumber
    \\
&= \max_{i} \left(
    \max_{x \in \cC_j}\ H^jA_ix \right) .
\end{align*}
We then have that,
through \eqref{eq:incr-linfty-hard-state}
and by rewriting $M_i = Q_i + k^*I$ for some $Q_i \geq 0$ and some $k^*\in \mathbb{R}$,
for each $i$
\begin{align*}
\max_{x \in \cC_j}\ H^j A_i x &= \max_{x \in \cC_j}\ (M_i H x)_j
= \max_{x \in \cC_j}\
(Q_iH x)_j + k^* H^j x \\
&\leq (Q_i \bm{1})_j + k^*
= (M_i \bm{1})_j = -\rate_w |B^T(H^j)^T|_1 .
\end{align*}

For \eqref{eq:linfty-state-output-final} we need to show that  
$  \|z\|_\infty= \|Cx\|_\infty \leq \rate_z $ for all trajectories $x$ such that $x(t) \in \mathcal{H}$.
Since  $\mathcal{H}$ is forward invariant, 
this is verified if $ [C; -C] x\leq  \rate_z \bm{1}$ for all $x$ such that  $Hx \leq \bm{1}$.
Indeed, from \eqref{eq:incr-linfty-hard-output},
$$
[C; -C] x = PH x \leq P \bm{1} = \eta_z \bm{1}.
\vspace*{-5mm}
$$
\hfill $\Box$

\begin{remark}
\cref{thm:incr-lone-gain,thm:incr-linfty-gain}
can be written in terms of the Jacobian $\partial f(x)$, without explicitly taking into account any convex relaxation $\mathcal{A}$.
This would make both theorems more general but not tractable, since we would have an infinite set of constraints to satisfy, one for each $x$.
In practice, to compute the gains, we would still need to rely on convex relaxations.
\end{remark}

\section{Computing Gain Bounds}
\label{sec:alg}
\subsection{Analysis}
\label{sec:anal}
Because of the connection between \eqref{eq:incr-lone-hard}
and \eqref{eq:incr-linfty-hard}, we only need one algorithm for computing both gain bounds.
We choose to design an algorithm for the $\Lone$ problem
and apply the algorithm to the adjoint system
if we are after $\Linfty$ gains bounds.
A comparison between our algorithm and other approaches is provided in \cref{sec:comparison}.

Because having variables multiplying each other makes a problem
nonconvex and nonlinear,
directly solving \eqref{eq:incr-lone-hard} with
the polyhedral candidate storage function represented by $V$
as a variable is not computationally tractable.
However,
for a fixed $V$,
\eqref{eq:incr-lone-hard}
can be efficiently solved using linear programming (LP).
This suggests a solution strategy based on an alternating iteration
to produce increasingly tight bounds on the gains
by locally optimizing a candidate $V$
(building on \cite{kousoulidis_polyhedral_2021}
that uses the same approach to find polyhedral Lyapunov functions
for autonomous systems).
We alternate between two steps, which we call
the \emph{gain estimation} step and the \emph{polyhedral modification} step.
Both steps are LP problems.
During the gain estimation step we consider $V$ fixed
and compute a bound on the gain.
During the polyhedral modification step,
we introduce a small variation in our candidate storage function
representation $\delta V$
with the goal of enabling a smaller gain bound to be verified.

\underline{Gain estimation step:} after fixing $V$,
variables $(\rate_w,P)$ and $(\rate_z,\cM)$ of \eqref{eq:incr-lone-hard}
become independent of each other.
We can then solve \eqref{eq:incr-lone-hard}
by solving two LP sub-problems: 
minimize $\rate_w$ subject to \eqref{eq:incr-lone-hard-input}
(with variables $\rate_w,P$),
and
maximize $\rate_z$ subject to \eqref{eq:incr-lone-hard-state}
(with variables $\rate_z,\cM$).
The optimal $\rate_w/\rate_z$ is then equal to
the solution to \eqref{eq:incr-lone-hard} for the given $V$.
Moreover,
if we fix $V$ such that 
$\psiV$ is a decaying incremental Lyapunov function for the internal dynamics,
\eqref{eq:incr-lone-hard-state} is guaranteed to be feasible
(since there must then exist $M_i$
such that $\offdiag(M_i) \geq 0, A_i V = V M_i$ and $1^TM_i < 0$)
and so the gain estimation step is guaranteed to obtain
a finite gain bound
(since \eqref{eq:incr-lone-hard-input} is always feasible
and implies $\rate_w > 0$).

\underline{Polyhedral modification step:}
we replace our objective $\min \rate_w/\rate_z$
 with $\min \log(\rate_w/\rate_z) = \min \{ \log(\rate_w)-\log(\rate_z) \}$
to obtain a simpler problem.
We then search for a small perturbation around 
solution
$(\rate_w, \rate_z, V, P, \cM)$
in the direction that most improves the linear approximation of our objective.
Working with small perturbations
allows us to express the modification step as a LP
in the (variational) variables
$
\delta \rate_w \in \R,\,
\delta \rate_z \in \R,\,
\delta V \in \R^{n \times m},\,
\delta P \in \R^{m \times 2n_w},\,
\delta \cM = \{\delta M_1 \in \R^{m \times m},\dots,
\delta M_k \in \R^{m \times m}\}
$
(where $k$ is the size of $\cA$)
and vector valued function of $V$,
$$\delta \hat{z}_V = [\sgn(CV_1)^T C\delta V_1,
\dots,
\sgn(CV_m)^T C\delta V_m]$$
(where $V_j$ and $\delta V_j$ are the $j^{th}$ columns of $V$ and $\delta V$ respectively):
\begin{subequations}
    \label{eq:big-LP}
    \begin{align}
       \min_{\delta \rate_w, \delta \rate_z,
       \delta V, \delta P, \delta \cM}
       \quad \> \delta \rate_w/\rate_w &- \delta \rate_z/\rate_z
       \label{eq:big-LP-objective}
       \\
        \noalign{\medskip}
       \begin{split}
        P + \delta P \geq 0,
        \     
        0 = \delta VP + V\delta P,
        \ 
        \delta \rate_w &\bm{1}^T = \bm{1}^T \delta P
       \end{split}
       \label{eq:big-LP-input}
            \\
        \noalign{\medskip}
        \begin{split}
        \offdiag(M_i + \delta M_i) \geq 0,
        \, 
        A_i \delta V = \delta V &M_i + V \delta M_i,
        \\
            \rate_z + \delta \rate_z > 0,
            \    
        -(\delta \rate_z \hat{z} + \rate_z \delta \hat{z}) &= \bm{1}^T \delta M_i,
        \\
        \text{for all } &i \in [k]
       \label{eq:big-LP-state}
        \end{split}
        \\
        \noalign{\medskip}
        |\vect(\delta V)|_\infty & \leq \varepsilon
        .
       \label{eq:big-LP-V}
    \end{align}
\end{subequations}
Breaking \eqref{eq:big-LP} into four parts,
\eqref{eq:big-LP-objective} is the linear approximation of
the variation in $\delta \log(\rate_w/\rate_z)$,
while \eqref{eq:big-LP-input} and \eqref{eq:big-LP-state} approximate
\eqref{eq:incr-lone-hard-input} and \eqref{eq:incr-lone-hard-state} respectively.
\cref{eq:big-LP-V} ensures that the variations considered are small
because it also bounds all other variables.

To account for the effect of the approximations
and ensure that the gain bound remains sound we only apply
the variation to the candidate storage function
to obtain $V_{\textrm{new}} := V + \delta V$
and then solve the gain estimation problem again for $V_{\textrm{new}}$.
For small enough $\varepsilon$,
this is guaranteed to lead to an improved gain bound.
We set $\varepsilon$ at
the $N^{th}$ iteration of our algorithm
to $\varepsilon_0/N$,
where $\varepsilon_0$ is
the initial step-size.
Additionally, if, after a polyhedral modification step,
a worse gain bound is computed
during the gain estimation step
or there exist redundant vertices in the polyhedral set
we halve $\varepsilon$ and
repeat the polyhedral modification step
(a vertex $V_j$ is redundant if
$\psiV(V_j) < 1$).
We stop our algorithm when $\varepsilon$
becomes smaller than a threshold.

Finally,
to bootstrap the gain estimation step to a feasible solution
we first use the algorithm from \cite{kousoulidis_polyhedral_2021}
to obtain a candidate $V$ that represents
a polyhedral incremental Lyapunov function for the internal dynamics.
This stage also controls the complexity
of Lyapunov functions, given by the number of vertices of $\cCV{V}$
(equivalently, by the number of columns of $V$). This number, $m$, is an input
to the algorithm.
In general, internal dynamics that are more oscillatory
will require functions with higher complexity
\cite{benvenuti_eigenvalue_2004}
(so if $m$ is set too low the algorithm in \cite{kousoulidis_polyhedral_2021}
might not be able to find a polyhedral Lyapunov function).
We also observe that higher complexity functions are better able to approximate
arbitrary gauge functions
and so can provide tighter bounds,
but this comes at an increased computational cost.
We briefly explore this trade-off in \cref{sec:motor-example}.

The full algorithms
for the $\Lone$ and $\Linfty$ incremental gains are summarized below.

\begin{algorithm}[htbp]
  {\bf Data:} Matrices $(\cA,B,C)$, \\
  polyhedral function complexity $\numrays$, \\
  and initial and minimum step-size parameters $(\varepsilon_0, \varepsilon_{min})$ \\
  {\bf Result:} Upper bound on $\Lone$ incremental gain, \\
  and corresponding storage function representation $V$ \\
  {\bf Procedure:} \\
  $(\varepsilon,N,\gamma',\delta V) := (\varepsilon_0, 1, \infty, 0)$ \\
  $V := V$ from \cite{kousoulidis_polyhedral_2021} LDI algorithm with inputs $(\cA,m)$ \\
  {\bf while }{\texttt{True}}: \\
  \hspace*{4mm} \emph{Gain Estimation Step:} \\
  \hspace*{4mm} $\rate_w, P :=
 \argmin_{\rate_w,P}\, \{\rate_w : \eqref{eq:incr-lone-hard-input}\}$
  using current $V$
  \\
  \hspace*{4mm} $\rate_z, \cM :=
 \argmax_{\rate_z,\cM}\, \{\rate_z : \eqref{eq:incr-lone-hard-state}\}$
  using current $V$
  \medskip \\
  \hspace*{4mm} {\bf if} $\rate_w/\rate_z > \gamma'$ or $\texttt{containsRedundant}(V)$: \\
  \hspace*{8mm} $\varepsilon,V := (\varepsilon/2, V-\delta V)$ \\
  \hspace*{4mm} {\bf else}: \\
  \hspace*{8mm} $(\varepsilon,N,\gamma') := (\varepsilon_0/N, N+1,\rate_w/\rate_z)$ \\
  \hspace*{4mm} {\bf if }{$\varepsilon < \varepsilon_{min}$}: \\
  \hspace*{8mm} return $(\gamma', V)$
  \medskip \\
  \hspace*{4mm} \emph{Polyhedral Modification Step:} \\
  \hspace*{4mm} $V := V + \delta V$ from 
  \eqref{eq:big-LP} using
    $(\rate_w, \rate_z, \varepsilon, V, P, \cM)$
  \caption{Incremental $\Lone$ gain bound computation \label{alg:lone}}
\end{algorithm}

\begin{algorithm}[htbp]
  {\bf Data:} Matrices $(\cA,B,C)$, \\
  polyhedral function complexity $\numrays$, \\
  and initial and minimum step-size parameters $(\varepsilon_0, \varepsilon_{min})$ \\
  {\bf Result:} Upper bound on $\Linfty$ gain, \\
  and corresponding invariant set representation $H$ \\
  {\bf Procedure:} \\
  $\gamma_\infty,H^T := $ \cref{alg:lone} with matrix inputs
  $(\cA^T, C^T, B^T)$ \\
  return $(\gamma_\infty, H)$
  \caption{Incremental $\Linfty$ gain bound computation \label{alg:linfty}}
\end{algorithm}

\subsection{Synthesis}
\label{sec:design}

After fixing $V$
the constraints in \eqref{eq:incr-lone-hard-state}
are convex with respect to the relaxed dynamics $\cA$.
This makes it easy to adapt the algorithms of
\cref{sec:anal} to linear feedback control synthesis
by adding extra free variables corresponding to the control parameters.

We consider
\eqref{eq:open-state-space}
with added
controlled inputs $u$
and
measured outputs $y$
\begin{align*}
    \dot{x} &= f(x) + B w + B_u u,\, 
    \\
    z &= C x,\\
    y &= C_y x .
\end{align*}
State feedback design is obtained by taking $C_y = I$ but our procedure 
works with any output matrix, allowing for the design of
static output feedback controllers of the form $u = Ky$.

The procedure is based on the extension 
of the gain estimation step of \cref{sec:anal}, by replacing
each constraint $A_i V = VM_i$ with
\begin{equation}
\label{eq:contr_syn_mod}
(A_i + B_u K C_y)V = VM_i .
\end{equation}
Since $K$ enters \eqref{eq:contr_syn_mod} linearly (for fixed $V$),
it is also straightforward to consider additional linear constraints on $K$
with the goal of enforcing bounds on the control gains or sparsity patterns.
By adding sums of constraints of the form
$\Gamma_l^T K \Gamma_r \leq g$, where $\Gamma_l$ and $\Gamma_r$ are matrices
of suitable size, 
we can limit the range of (linear combinations of) the coefficients of $K$.  
Likewise, constraints of the form
$\Gamma_l^T K \Gamma_r = 0$
can be used to enforce specific sparsity patterns on $K$. 

The introduction of the additional variable $K$ also imposes minor changes to
the polyhedral modification step.
The small perturbations in \eqref{eq:big-LP} must
take into account the extended solution manifold $$(\rate_w, \rate_z, V, P, \cM, K) .$$
This means that we need to add $\delta K$ as a new variable to \eqref{eq:big-LP}
and replace the constraints
$A_i \delta V = \delta V M_i + V\delta M_i$
with
\begin{equation}
(A_i + B_u K C_y)\delta V + B_u \delta K C_y V
= \delta V M_i + V\delta M_i .
\end{equation}
We also add constraint $\|\vect(\delta K)\|_\infty \leq \varepsilon$
to guarantee that only small perturbations are considered.
Finally, we add constraints
$\Gamma_l^T (K+\delta K) \Gamma_r \leq g$
and $\Gamma_l^T (K+\delta K) \Gamma_r = 0$
for each constraint
$\Gamma_l^T K \Gamma_r \leq g$ and $\Gamma_l^T K \Gamma_r = 0$,
respectively.

\subsection{Comparison with Other Approaches}
\label{sec:comparison}
The use of polyhedral sets for the computation of gains
is not new. They have been used in
\cite[Section 6.4]{blanchini_set-theoretic_2015}
and \cite{fialho_l1_1995} for the computation of 
the $\Linfty$ gain of discrete time LTV/LPV systems.

In \cite{blanchini_set-theoretic_2015},
the algorithm computes the largest set included in $\|Cx\|_\infty \leq \mu$
that is invariant under all possible $\|u\|_\infty \leq 1$.
This is done by starting with polyhedral set
$\|Cx\|_\infty \leq \mu$
and propagating the set backwards in time until the set becomes invariant,
in which case $\mu$ is an upper bound for the $\Linfty$ gain,
or no longer contains $0$,
in which case $\mu$ is a lower bound for the $\Linfty$ gain.
This is repeated for different values of $\mu$ to obtain a sharp bound.

The approach of \cite{fialho_l1_1995} 
computes the reachable set from the origin by forward propagation
in time.
After each step, the smallest $\mu$ for which set
$\|Cx\|_\infty \leq \mu$ contains the reachable set
is a lower bound on the gain.
\cite{fialho_l1_1995} also provide bounds on the error
between the computed lower bound and the actual gain of the system,
as a function of the number of forward steps taken and of the
system spectral radius. This allows the algorithm
to terminate in a pre-determined finite number of steps
as a function of the required accuracy.
A drawback of the approach is that it scales worse
for inputs of large dimension.
This is because V-rep polyhedra must be used throughout,
including to represent the $\Linfty$ unit ball
of the image of the input $Bw$,
which by itself requires $2^{n_w}$ vertices.

The approaches of both \cite{blanchini_set-theoretic_2015} and \cite{fialho_l1_1995}
can be applied to continuous time systems
through
a forward Euler discretization scheme
(from $\dot{x} = Ax$ to $x^+ = [I + \tau A]x$).
Discretizing the systems in this way
guarantees that the set invariance properties
necessary to bound the $\Linfty$ gain
will hold for the continuous system if they hold for the discretized system
\cite{blanchini_set-theoretic_2015}.
By using smaller discretization time constants the bound on the gain
can then be made increasingly tighter, but
at an increasing computational cost.
This is a sharp difference with our approach, which is based on conditions
that apply directly to continuous time systems, without the need of any discretization step. 

The approaches of \cite{blanchini_set-theoretic_2015}
and \cite{fialho_l1_1995}
iteratively shrink or grow polyhedra at every step of computation. This may lead to
increasingly large matrices to represent the polyhedra,
with the matrix size growing exponentially with the number of steps taken.
In contrast, our approach does not shrink or grow polyhedra but reshapes them at each step,
keeping constant the size of the representation matrix.
This is important since the computational cost at each
iteration scales with the size of the polyhedral representation.
As such,
while for the approaches of \cite{blanchini_set-theoretic_2015}
and \cite{fialho_l1_1995}
there is guaranteed convergence to the actual LTV/LPV gain in a finite number of steps
for discrete or discretized systems,
each step could become computationally intense, with
no upper bound obtained until the process terminates.
In contrast, there is no guaranteed convergence to the actual LTV/LPV gain
for our approach but 
each step always has a fixed cost and produces a
tighter upper bound. The complexity of the polyhedra, characterized by
the number of vertices $m$, is governed by the user.

Another important unique strength of our optimization-based approach
is that it allows for synthesis problems to be tackled directly.

\section{Examples}
\label{sec:examples}

\appendtographicspath{{examples/motor_example/code+figs/}}
\subsection{DC Motor}
\label{sec:motor-example}
We revisit the DC motor example of  \cite{kousoulidis_polyhedral_2021}
for the computation of incremental system gains.
In the first part, we estimate a bound on the incremental $\Linfty$ gain between
a disturbance acting on an open loop (linear) DC motor and its output speed. Later, 
we design a feedback controller that stabilizes the motor position while 
minimizing the effect of disturbances in closed loop.

The nominal parameters are:
inertia $J_0=.01$, viscous friction constant $b_0=.1$,
electromagnetic field constant $E_0=.01$, resistance $R_0=1$, and inductance $L_0=.5$.
In each case,
we run the algorithm 10 times with
$\epsilon_0 = 1/5$ and 
different random seeds for $V$
and report the result of the best run.

\begin{figure}[tbp]
    \centering
   \begin{minipage}{.49\linewidth}
        \includegraphics[width=\linewidth]{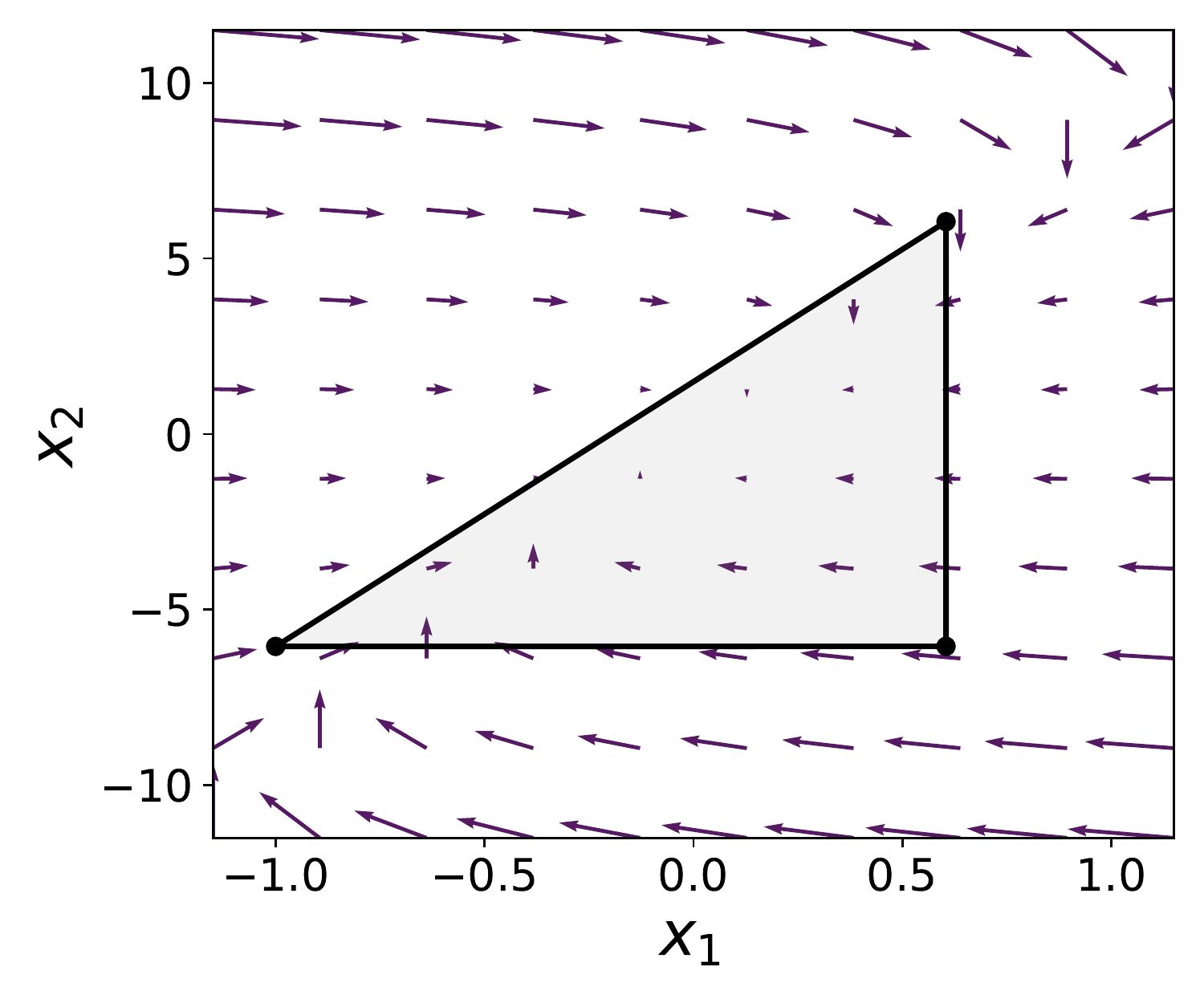}
   \end{minipage}
   \begin{minipage}{.49\linewidth}
        \includegraphics[width=\linewidth]{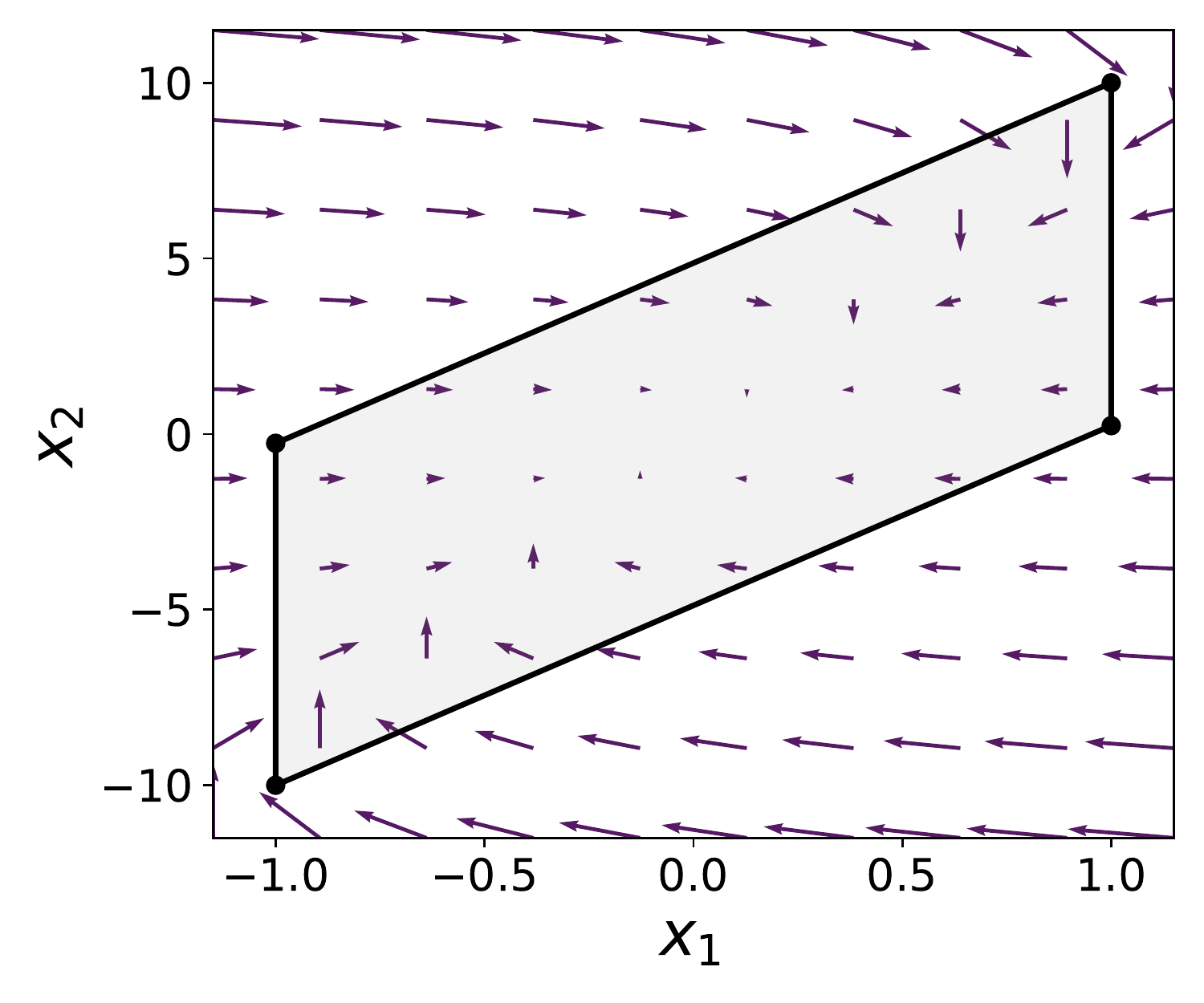}
   \end{minipage}
    \caption{
    Three and four sided H-rep sets
    found using \cref{alg:linfty}
    bounding the $\Linfty$ gain of \eqref{eq:motor_speed}
    (with nominal parameter values)
    to $\leq 0.083$ and $\leq 0.050$ respectively.
    }
   \label{fig:motor_speed_nominal}
 \end{figure}

For the analysis scenario, the variational dynamics of the DC motor are
\begin{align}
    \dot{\delta x}
    &= \begin{bmatrix}
        -b/J & E/J \\
        -E/L & -R/L
    \end{bmatrix}
    \delta x
    +
    \begin{bmatrix}
        0 \\
        1
    \end{bmatrix}
    \delta w
    \label{eq:motor_speed} \\
    \delta z &=
    \begin{bmatrix}
        1 & 0
    \end{bmatrix}
    \delta x \nonumber
    ,
\end{align}
where the first state corresponds to the speed of the motor and the second to the current.
We first consider the case where each parameter in \eqref{eq:motor_speed} is set to its nominal value,
$b = b_0$, $E = E_0$, $J = J_0$, and $R = R_0$. 
We apply
\cref{alg:linfty}
to find H-rep sets with three and four half-spaces
(shown in \cref{fig:motor_speed_nominal})
for which the bounds obtained from \cref{thm:incr-linfty-gain}
are around $0.083$ and $0.050$ respectively.
In this case, because the impulse response is positive,
the actual $\Linfty$ gain can be analytically computed to be $1/20.02$ 
using \cref{remark:pos}
and matches the bound for the four half-space set
to 4 significant figures.

\begin{figure}[tbp]
    \centering
   \begin{minipage}{.49\linewidth}
        \includegraphics[width=\linewidth]{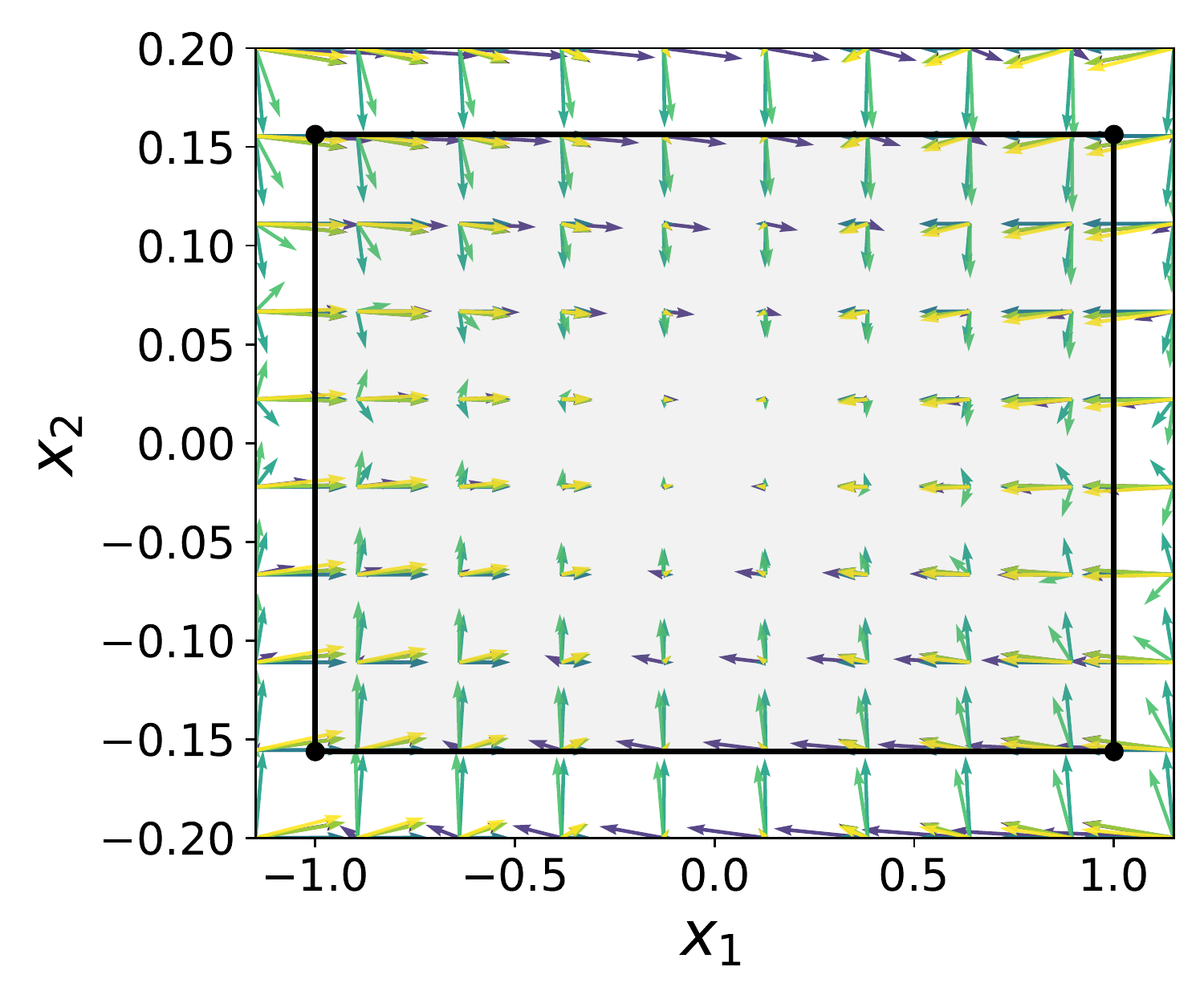}
   \end{minipage}
   \begin{minipage}{.49\linewidth}
        \includegraphics[width=\linewidth]{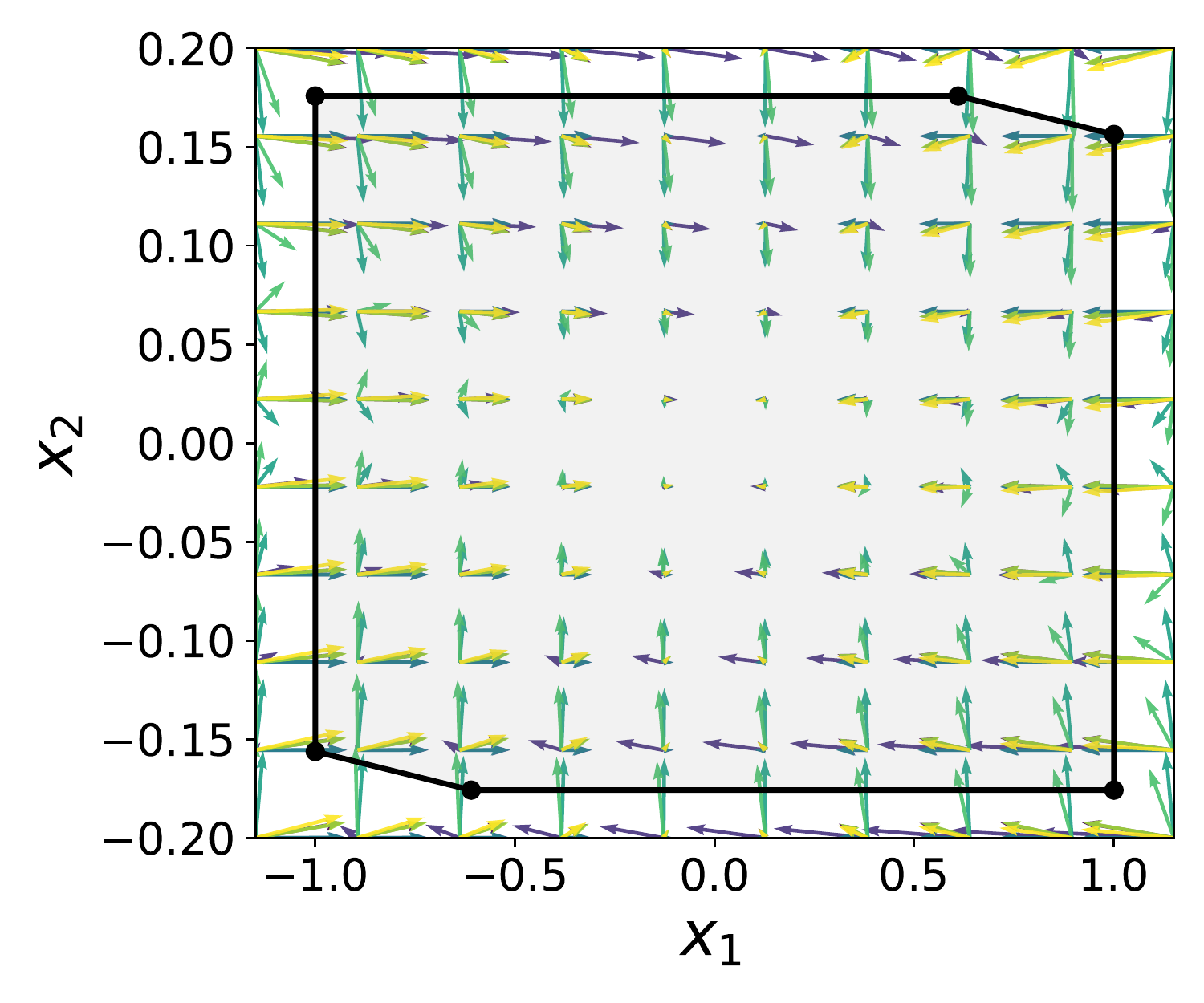}
   \end{minipage}
   \begin{minipage}{.49\linewidth}
        \includegraphics[width=\linewidth]{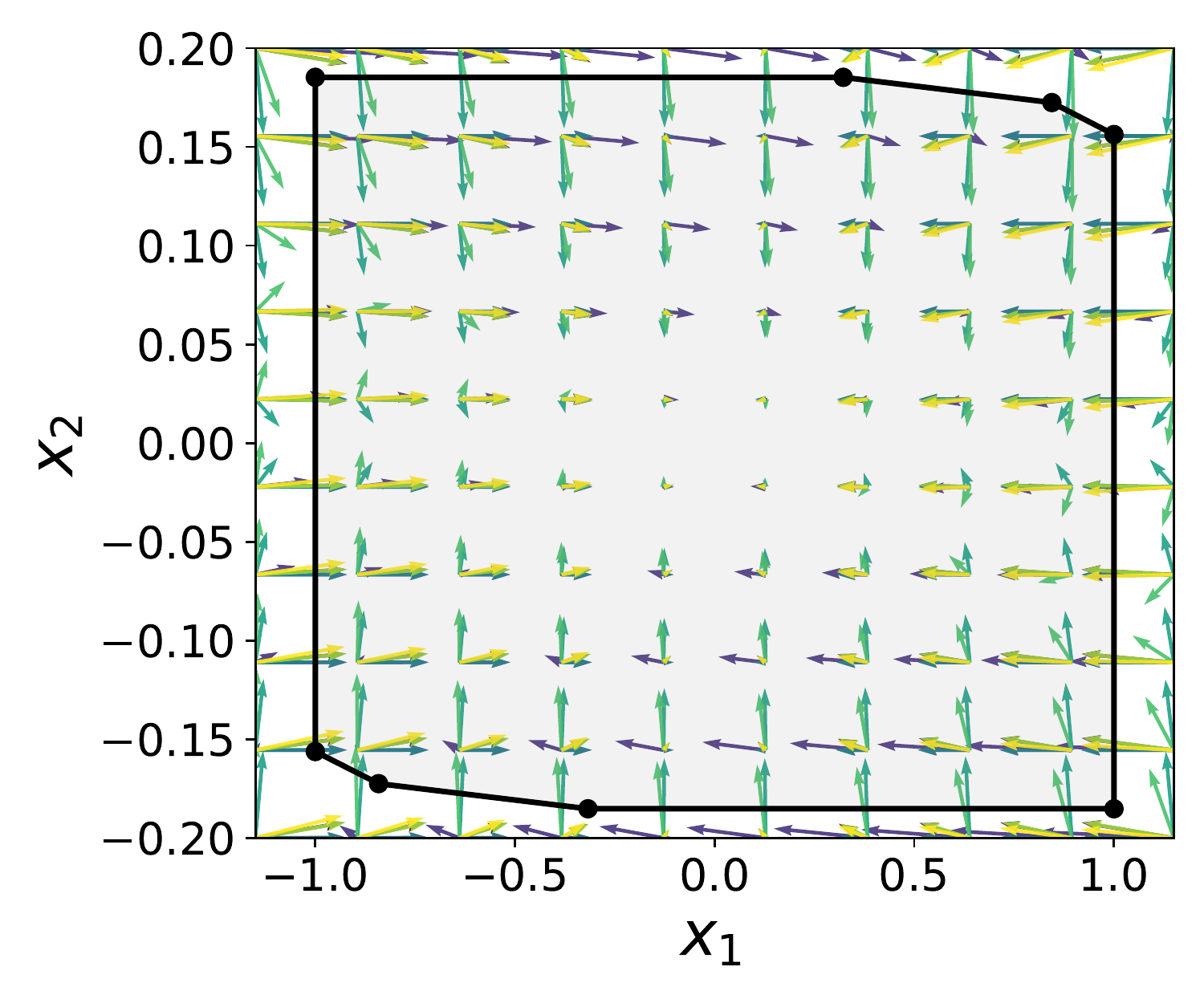}
   \end{minipage}
   \begin{minipage}{.49\linewidth}
        \includegraphics[width=\linewidth]{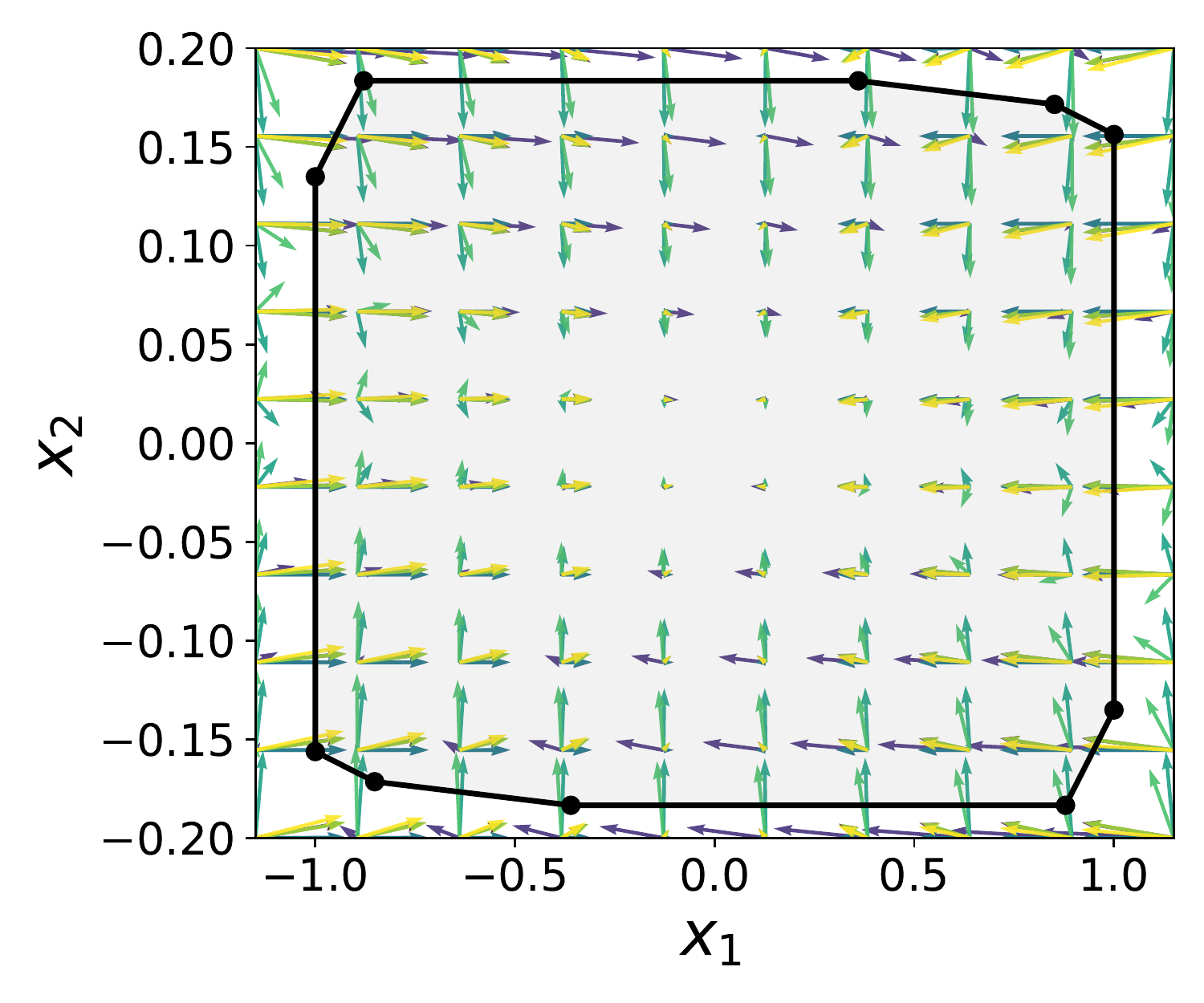}
   \end{minipage}
    \caption{
    Four, six, eight, and ten sided
    H-rep polyhedral sets
    found using \cref{alg:linfty}
    bounding the incremental $\Linfty$ gain of
    \eqref{eq:motor_speed}-\eqref{eq:motor_speed_params}
    (with $\delta_s = 8$) to
    $6.6$, $5.2$, $4.8$, and $4.4$ respectively.
    }
   \label{fig:motor_speed_robust}
 \end{figure}

We then consider \eqref{eq:motor_speed}
with uncertain parameters
\begin{equation}
J\! \in\! [J_0/\delta_s,\delta_s J_0],\, b\! \in\! [b_0/\delta_s,\delta_s b_0],\, E\! \in\! [E_0/\delta_s,\delta_s E_0] 
\label{eq:motor_speed_params}
\end{equation}
for $\delta_s = 8$,
leading to a relaxation set $\cA$ composed of $8$ vertices.
Because this is not a linear system,
the incremental $\Linfty$ gain can't be computed analytically in this case.
A lower bound for it can be obtained
by computing the gain for linear systems
with fixed parameters picked from the vertices of the sets in \eqref{eq:motor_speed_params}.
Out of these,
the highest occurs when $(J,b,E) = (J_0/\delta_s,b_0/\delta_s,\delta_s E_0)$
and can be analytically computed to be around $2.12$.
We use \cref{alg:linfty} to compute an upper bound.
We work with H-rep sets
with four, six, eight, and ten sides,
leading to upper bounds on the incremental gain of
$6.6$, $5.2$, $4.8$, and $4.4$ respectively.
This shows how increasing the complexity of the polyhedral representation leads to
less conservative results.
The polyhedra found are shown in \cref{fig:motor_speed_robust}.

For the synthesis scenario we design a linear state feedback controller and a linear output feedback controller, the latter
using motor position and motor current. The goal is to minimize the worse case effect that bounded disturbance forces 
have on the motor position. 

\begin{figure}[tbp]
   \begin{minipage}{.49\linewidth}
        \includegraphics[width=\linewidth]{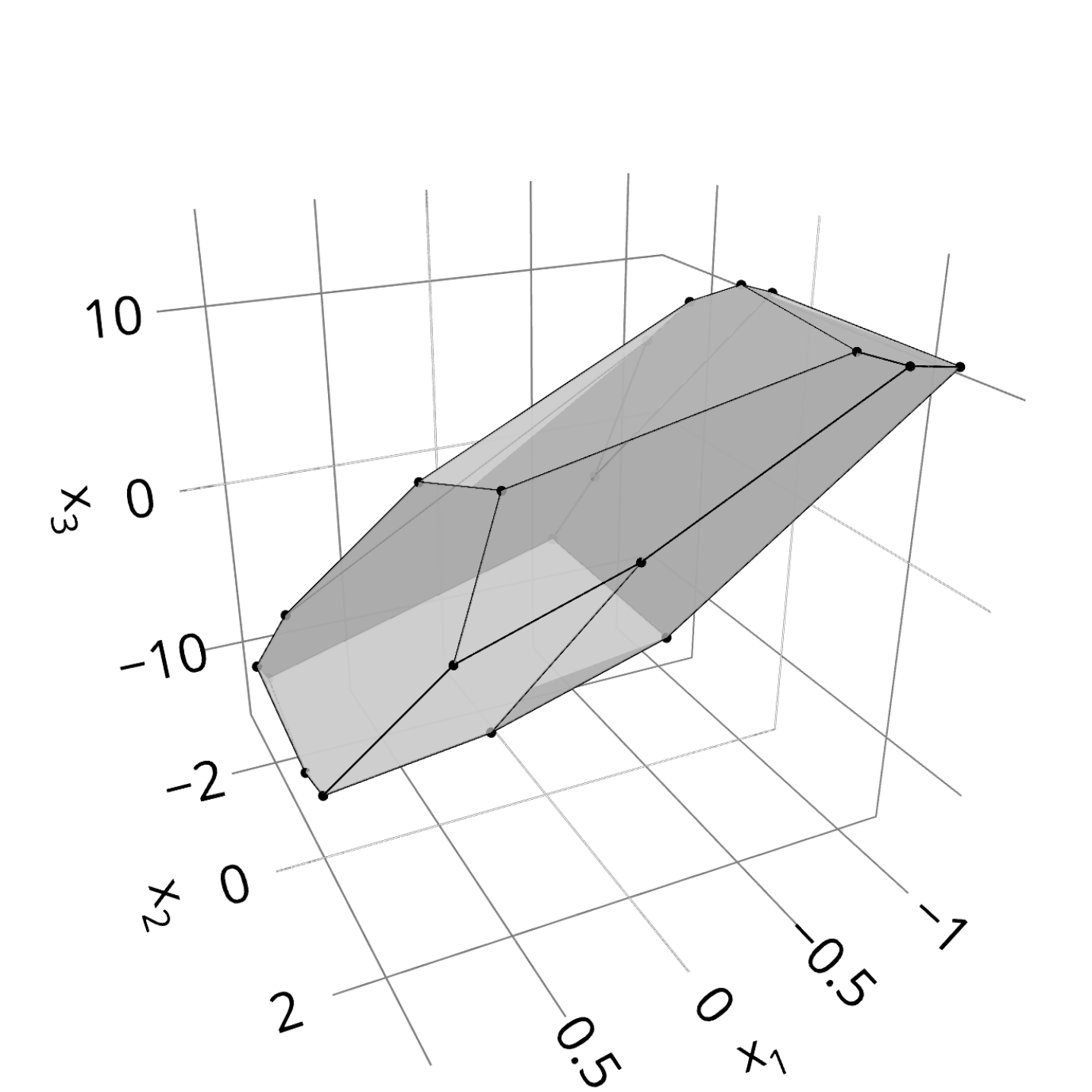}
   \end{minipage}
   \begin{minipage}{.49\linewidth}
        \includegraphics[width=\linewidth]{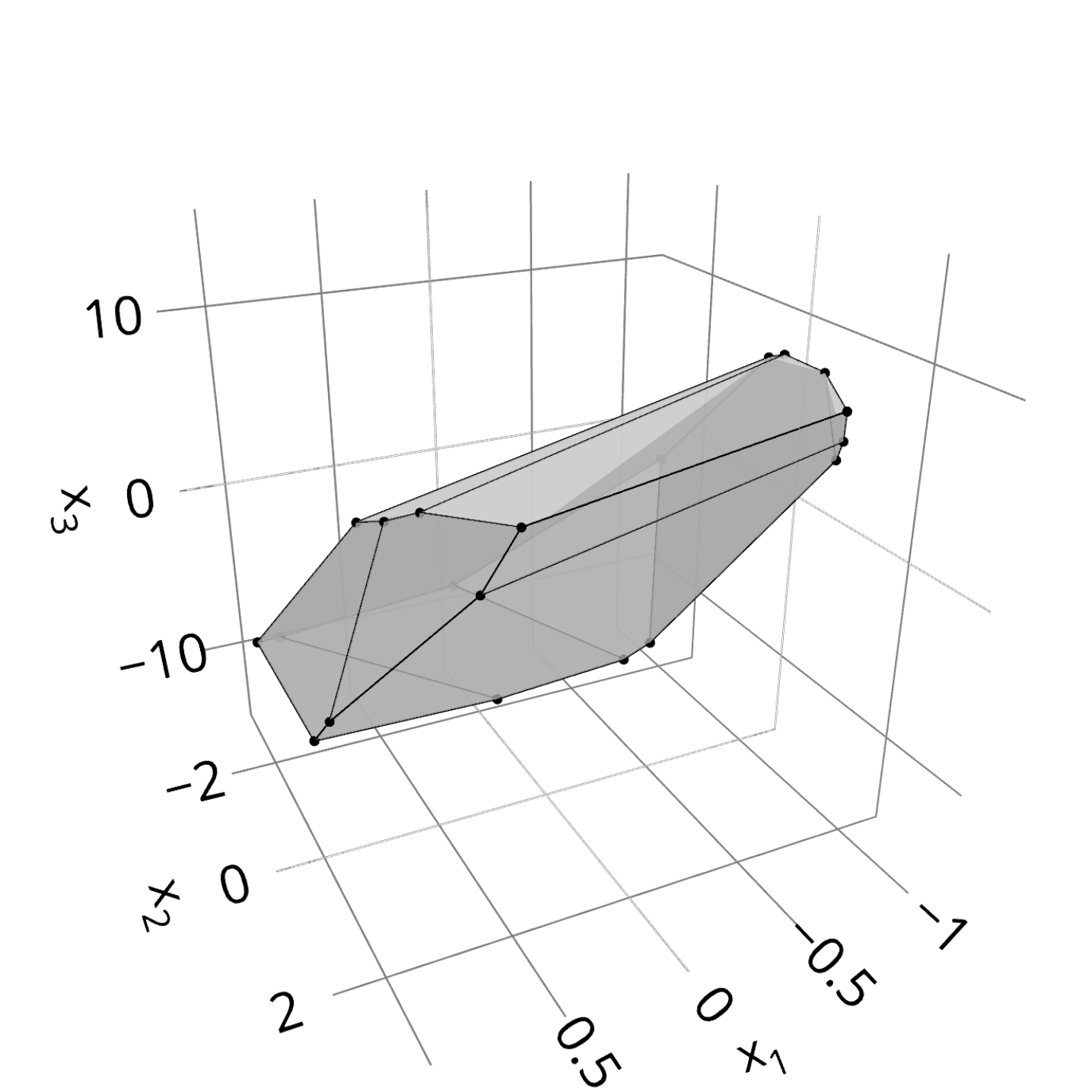}
   \end{minipage}
   \caption{
    12 sided H-rep sets found 
    simultaneously to feedback schemes
    \eqref{eq:motor_position_full},\eqref{eq:motor_position_feedback_full} (left)
    and
    \eqref{eq:motor_position_no_speed},\eqref{eq:motor_position_feedback_no_speed} (right)
    using \cref{alg:linfty} adapted to synthesis.
    They bound the incremental $\Linfty$ gain for 
    \eqref{eq:motor_position}-\eqref{eq:motor_position_params} (with $\delta_p = 1.4$)
    to $\leq 0.17$
    and
    $\leq 0.26$
    respectively.
   }
   \label{fig:motor_position}
 \end{figure}

The variational dynamics reads 
\begin{align}
    \dot{\delta x}
    &= \begin{bmatrix}
        0 & 1 & 0 \\
        0 & -b\!/J & E\!/J \\
        0 & -E\!/L & -R\!/L
    \end{bmatrix}
    \delta x
    + \!\begin{bmatrix}
        0 \\
        0 \\
        1\!/L
    \end{bmatrix}
    \delta u
    + \!\begin{bmatrix}
        0 \\
        1 \\
        0
    \end{bmatrix}
    \delta w
    \label{eq:motor_position} \\
    \delta z &=
    \begin{bmatrix}
        1 & 0 & 0
    \end{bmatrix}
    \delta x
    ,
    \label{eq:motor_position_output}
\end{align}
where the first and second states are related to the motor position and velocities, and the
last state is related to the motor current. We also consider the following parametric uncertainties
\begin{equation}
J \in [J_0/\delta_p,\delta_pJ_0],\, E \in [E_0/\delta_p,\delta_pE_0]
\label{eq:motor_position_params}
\end{equation}
where $\delta_p = 1.4$.
We first look into state feedback, using the output
\begin{equation}
y_{s} = - I x .
\label{eq:motor_position_full}
\end{equation}
We then move to output feedback design, from
\begin{equation}
y_{o} = - \begin{bmatrix}
    1 & 0 & 0 \\
    0 & 0 & 1
    \end{bmatrix} x .
\label{eq:motor_position_no_speed}
\end{equation}
In both case, we adopt a polyhedral representation given by
12 half-spaces. With the modified conditions of Section \ref{sec:design},
our algorithm proves a closed-loop $\Linfty$ bound of $0.17$, for
the state feedback controller 
\begin{equation}
 u = 
\begin{bmatrix}
    27.5 & 3.04 & 1.41
\end{bmatrix}
y_{s} .
\label{eq:motor_position_feedback_full}
\end{equation}
Likewise, the algorithm derives a closed-loop $\Linfty$ gain of $0.26$ 
for the output feedback controller
\begin{equation}
 u = 
\begin{bmatrix}
    21.8 & 1.74
\end{bmatrix}
y_{o} .
\label{eq:motor_position_feedback_no_speed}
\end{equation}
The associated polyhedra are shown in \cref{fig:motor_position}.

\appendtographicspath{{examples/compartmental_system/code+figs/}}
\subsection{Compartmental System}
We next consider a small compartmental system example.
Compartmental systems are naturally suited for pharmacological and epidemiological models,
water reservoirs and queuing systems,
and as discretizations of partial differential equations (PDEs)
\cite{haddad_nonnegative_2010}.
Each compartment can have inflows and outflows, both internal and external.
In our example we consider a basic compartmental system where 
matter moves among compartments and its quantity is conserved
through internal flows but can be added or removed from the overall system through external inflows or outflows.
We will once again be running the algorithm 10 times
with $\varepsilon_0 = 1/5$ and different random seeds for $V$
and reporting the result of the best run.

Specifically, we consider the design of output feedback controllers with integral action
for a compartmental model that
approximates diffusion through a pipe.
The pipe is modelled by 5 compartments. The size of the state space of the closed-loop system is $6$, 
with the sixth state needed for the integral action of the controller. The objective is to minimize the effects
of external disturbances on the measured output associated to the last compartment. 
We express this as incremental $\Lone$ and $\Linfty$
gain synthesis problems for system
\begin{equation}
\begin{split}
\dot{x}_1 &= - \sigma x_1 + x_2 + w_1 + u \\
\dot{x}_2 &= \sigma x_1 -2x_2 + x_3 + w_2 \\
\dot{x}_3 &= x_2 -2x_3 + x_4 + w_3 \\
\dot{x}_4 &= x_3 -2x_4 + x_5 + w_4 \\
\dot{x}_5 &= x_4 - 2x_5 + w_5 \\
\dot{q} &= r-x_5 \\
z & =  x_5 - r ,
\end{split}
\end{equation}
where each state $x_i$ corresponds to the matter in the $i$th compartment,
and $q$ corresponds to the integral error $r-x_5$, for the desired reference value $r$.
The input $u$ represents the rate at which matter is being injected at the beginning of the pipe.
We also set
$
\sigma \in [0.5,1]
$,
encompassing dynamics where the flow
from the first compartment to the second becomes partially obstructed.
Matter is conserved everywhere except in the last compartment,
where it is escaping with a rate of $1$. Variables $w_i$ represent independent external disturbances
that are acting on each compartment.
The proportional-integral controllers $u = K y$ use the  measured output
\begin{equation}
y = 
\begin{bmatrix}
-x_1 & 
r-x_5 &
q
\end{bmatrix}^T
.
\label{eq:compartmental-observations}
\end{equation}
The variational dynamics associated to the system reads 
\begin{equation}
    \label{eq:compartmental_main}
\begin{split}
\dot{\delta x}_1 &= - \sigma \delta x_1 + \delta x_2 + \delta w_1 + \delta u \\
\dot{\delta x}_2 &= \sigma \delta x_1 -2\delta x_2 + \delta x_3 + \delta w_2 \\
\dot{\delta x}_3 &= \delta x_2 -2\delta x_3 + \delta x_4 + \delta w_3 \\
\dot{\delta x}_4 &= \delta x_3 -2 \delta x_4 + \delta x_5 + \delta w_4 \\
\dot{\delta x}_5 &= \delta x_4 - 2\delta x_5 + \delta w_5 \\
\dot{\delta q} &= -\delta x_5 \\
\delta z & = \delta x_5 ,
\end{split}
\end{equation}
which, for 
$
\sigma \in [0.5,1]
$,
leads to a relaxation set $\cA$ composed of 2 vertices.

We use \cref{alg:lone,alg:linfty} 
to design controllers that minimize the 
incremental $\Lone$ and $\Linfty$ gains from the disturbance input $w$ to the performance output $z$.
The results are summarized in \cref{tbl:compartmental_robust_gains}.

\begin{table}[h]
    \centering
    \begin{tabular}{c | c | c | c}
        & $K$
        & {\small Upper Bound}
        & {\small Lower Bound}
        \\
        \hline
        $\Lone$
        &
        $\begin{bmatrix}
            2.51 & 1.28 & 0.27
        \end{bmatrix}
        $
        & 2.47
        & 1.47
        \\
        $\Linfty$
        &
        $\begin{bmatrix}
            9.10 & 6.45 & 0.23
        \end{bmatrix}
        $
        & 5.06
        & 3.54
    \end{tabular}
    \caption{Table of control parameters
    and incremental gain bounds found using \cref{alg:lone,alg:linfty}
    (with $m = 18$)
    for minimizing each target gain bound for the linear case.
    The lower bounds are computed from
    the vertices of the closed-loop linear systems
    satisfying \eqref{eq:compartmental_main}
    with $\sigma \in [.5,1]$ using \eqref{eq:lone-tf} and \eqref{eq:linfty-tf} combined with \cite{rutland_computing_1995}.
    }
    \label{tbl:compartmental_robust_gains}
\end{table}

Both performance specifications have interesting applications.
For example, 
the $\Linfty$ gain bound bounds the maximum amount by which
the regulated output will vary about its reference value
when subjected to independent bounded disturbances
on each of the compartments.

Beyond their direct interpretations in terms of performances, 
the gains also quantify the robustness of the closed-loop system
to unmodeled dynamics, via the
small gain theorem. For example, the gain in Table \ref{tbl:compartmental_robust_gains}
guarantees closed-loop incremental stability
even if the last compartment is leaking at a rate
$\in [0.6,1.4]$. This is obtained by setting $\delta w_5$ to $\pm (1/\gamma) \delta x_5$,
where $\gamma$ is the smallest upper bound in \cref{tbl:compartmental_robust_gains}.

\section{Conclusions}
\label{sec:conclusions}
We have discussed theoretical conditions for the estimation of bounds on the incremental $\Lone$ and $\Linfty$
gains of nonlinear systems. Their computation is performed via a two-step procedure based on linear programming problems.
To the best of our knowledge, the characterization of $\Lone$ gain bounds is a novel contribution 
while the $\Linfty$ bounds are closer to results available in the literature, based on set invariance.
The algorithms outlined in the paper are sound but susceptible to local optimal. 
They can be used both for system analysis and for closed-loop control synthesis.
The gains computed can quantify performance specifications
or be used in conjunction with the small gain theorem
to prove stability of interconnections. The practical use of the
algorithms was illustrated by two numerical examples. 

More work is needed on optimizing, tweaking, and formally evaluating
the properties of the optimization procedure used in the algorithm.
From the modes of failure observed
in sub-optimal solutions, it seems that explicitly enforcing symmetry
and adding the dual constraints of
the solutions from the gain estimation step
to the polyhedral modification step helps avoid a big class of local optima.
Other extensions include weighted norms and synthesis of nonlinear controllers \cite{blanchini_persistent_1995}.
This is particularly relevant because it is well known that linear state feedback controllers are sub-optimal for the $\Linfty$ problem \cite{diaz-bobillo_state_1992}.

\bibliographystyle{abbrv}
\bibliography{22_poly_open_conf}
                                                   
\end{document}